\documentclass[12pt]{article}
\usepackage{amsmath,amssymb,amsthm, amsfonts}
\usepackage{enumerate}
\usepackage{hyperref}
\usepackage{graphicx}
\usepackage{url}
\usepackage[mathlines]{lineno}
\usepackage{dsfont} 
\usepackage{soul}
\usepackage{tikz, graphicx, subcaption, caption}
\usetikzlibrary{positioning}
\usepackage{color}
\definecolor{red}{rgb}{1,0,0}

\definecolor{blue}{rgb}{0,0,1}

\definecolor{green}{rgb}{0,.6,0}

\newtheorem{thm}{Theorem}[section]
\newtheorem{cor}[thm]{Corollary}
\newtheorem{lem}[thm]{Lemma}
\newtheorem{prop}[thm]{Proposition}

\newtheorem{obs}[thm]{Observation}

\newtheorem{quest}[thm]{Question}

\theoremstyle{definition}
\newtheorem{rem}[thm]{Remark}

\theoremstyle{definition}
\newtheorem{defn}[thm]{Definition}

\theoremstyle{definition}
\newtheorem{ex}[thm]{Example}

\numberwithin{figure}{section}   
\numberwithin{table}{section}   
\numberwithin{equation}{section}

\newcommand{\ZZ}{\mathbb{Z}}

\newcommand{\G}{\mathcal{G}}

\newcommand{\Z}{\operatorname{Z}}
\newcommand{\M}{\operatorname{M}}
\newcommand{\ZIR}{\operatorname{ZIR}}
\newcommand{\zir}{\operatorname{zir}}
\newcommand{\pd}{\gamma_P}
\newcommand{\tw}{\operatorname{tw}}
\newcommand{\zbar}{\ol\Z}
\newcommand{\fr}{\operatorname{Fr}}
\newcommand{\IR}{\operatorname{IR}}
\newcommand{\ir}{\operatorname{ir}}

\newcommand{\nul}{\operatorname{null}}
\newcommand{\Gc}{\overline{G}}
\newcommand{\dist}{\operatorname{dist}}

\newcommand{\noi}{\noindent}

\newcommand{\bit}{\begin{itemize}}
\newcommand{\eit}{\end{itemize}}
\newcommand{\ben}{\begin{enumerate}}
\newcommand{\een}{\end{enumerate}}
\newcommand{\beq}{\begin{equation}}
\newcommand{\eeq}{\end{equation}}
\newcommand{\bea}{\begin{eqnarray*}}
\newcommand{\eea}{\end{eqnarray*}}
\newcommand{\bean}{\begin{eqnarray}}
\newcommand{\eean}{\end{eqnarray}}
\newcommand{\bpf}{\begin{proof}}
\newcommand{\epf}{\end{proof}\ms}
\newcommand{\bmt}{\begin{bmatrix}}
\newcommand{\emt}{\end{bmatrix}}
\newcommand{\ms}{\medskip}

\newcommand{\beqa}{\begin{array}}
\newcommand{\eeqa}{\end{array}}
\newcommand{\OL}{\overline}
\newcommand{\lc}{\left\lceil}
\newcommand{\rc}{\right\rceil}
\newcommand{\lf}{\left\lfloor}
\newcommand{\rf}{\right\rfloor}
\newcommand{\lp}{\left(}
\newcommand{\rp}{\right)}

\newcommand{\wh}{\widehat}

\newcommand{\du}{\mathbin{\,\sqcup\,}}
\newcommand{\ol}{\overline}

\title{Zero forcing irredundant sets}
\author{Bryan A.\ Curtis\thanks{Department of Mathematics, Iowa State University,
Ames, IA 50011, USA (dr.bryan.curtis.math@gmail.com)}\and  Leslie Hogben\thanks{American Institute of Mathematics, Caltech 8-32, 1200 E California Blvd., Pasadena, 91125 USA
(hogben@aimath.org) and Department of Mathematics, Iowa State University,
Ames, IA 50011, USA.} \and Adriana Roux\thanks{Corresponding author; Department of Mathematical Sciences, Stellenbosch University, Stellenbosch, 7600, South Africa (rianaroux@sun.ac.za)}}

\begin{document}
\maketitle
\vspace{-8pt}

\begin{abstract}  Irredundance has been studied in the context of dominating sets, via the concept of a private neighbor.  Here irredundance of zero forcing sets is introduced via the concept of a private fort and the upper and lower zero forcing irrdedundance numbers $\ZIR(G)$ and $\zir(G)$ are defined. Bounds on  $\ZIR(G)$ and $\zir(G)$ are established and graphs having extreme  values of $\ZIR(G)$ and $\zir(G)$ are characterized. The effect of the join and corona operations is studied. As the concept of a zero forcing irrdedundant set is new, there are many questions for future research.
\end{abstract}

\noi {\bf Keywords} zero forcing; irredundance; ZIr-set; ZIR

\noi{\bf AMS subject classification}  05C69, 05C50, 05C57


\section{Introduction}\label{s:intro} 

In this paper we introduce an analog of irredundant sets for zero forcing, called zero forcing irredundant sets, with a private fort playing the role of a private neighbor.  We also investigate the relationships between the zero forcing and zero forcing irredundance numbers in analogy with  the Domination Chain (see \eqref{eq:domchain}).  We begin this introduction with a brief discussion of domination and irredundance to motivate our work on zero forcing irredundant sets, followed by definitions of zero forcing and forts.  We then introduce the definition of zero forcing irredundance.  Finally we outline the structure of the rest of the paper and summarize notation we will use.

\subsection{Motivation}
A set $S$ is a \emph{dominating  set} of a graph $G$ if every vertex of $G$ is in $S$ or a neighbor of a vertex in $S$.  A set $S$ of vertices in a graph $G$ is \emph{irredundant} if each vertex in $S$ dominates a vertex of $G$ that is not dominated by any other vertex in $S$  (this may be a neighbor or itself), called a \emph{private neighbor}. The study of irredundance was introduced in 1978 by Cockayne, Hedetniemi and Miller \cite{CHM78} in relation to the study of minimal dominating sets. Irredundance has been studied extensively since then (see \cite{MR21} and the references therein). 

The \emph{domination number} $\gamma(G)$ and the \emph{upper domination number} $\Gamma(G)$ of $G$ are the minimum and maximum cardinalities of a minimal dominating set of $G$. The \emph{upper irredundance number} $\IR(G)$ and the \emph{lower irredundance number} $\ir(G)$ of $G$ are the maximum and minimum cardinalities of a maximal irredundant set of $G$.
A set that is both dominating and irredundant is a minimal dominating and maximal irredundant set. In fact, minimal dominating sets have been characterized as exactly those sets that are both dominating and irredundant \cite{CHM78}.

The relationships between the upper and lower parameters for domination and irredundance  are expressed in the  \emph{Domination Chain}  introduced in \cite{CHM78}:   
\beq \label{eq:domchain}\ir{G} \leq \gamma(G) \leq \Gamma(G) \leq \IR(G).\eeq
 
Various aspects of the Domination Chain have been studied, including complexity of computing parameters in the Domination Chain (see, for example, \cite{BBCF20domchain} and the references therein). The conditions under which these parameters are equal have been studied extensively (see \cite{MR21} for a survey).

\subsection{Zero forcing}
Zero forcing was introduced as a bound for the maximum nullity among real symmetric matrices whose pattern of nonzero off-diagonal entries is described by the edges of a graph in \cite{AIM} and independently in other applications.

Zero forcing is a propagation process on a graph $G$.  Starting with an initial set of blue vertices, the process colors vertices blue by repeated applications of the \emph{color change rule}: A blue vertex $u$ can change the color of a white vertex $w$ to blue if $w$ is the only white neighbor of $u$. A \emph{zero forcing set}   of  $G$  is a subset of vertices $B$ such that  if $B$ is the initial set of blue vertices, then the repeated application of the color change rule will cause the whole graph to turn blue.
The \emph{zero forcing number} of  $G$ is  the minimum cardinality of a zero forcing set.  
The  \emph{upper zero forcing number} $\zbar(G)$ of a graph $G$ is the maximum cardinality of a minimal zero forcing set of $G$.

\subsection{Private forts}

Let $G$ be a  graph. A nonempty set $F\subseteq V(G)$ is a \emph{fort} if for all $v\in V(G)\setminus F$, $|F\cap N(v)|\ne 1$ where $N(v)$ denotes  the set of neighbors of $v$.  Note that $V(G)$ is a fort (vacuously) and in a connected graph of order at least two, a fort must contain at least two vertices.
Forts present obstructions to zero forcing: 
Given a set $B$ of blue vertices, some vertex in $B$  can perform a force in $G$ if and only if $V(G)\setminus B$ is not a fort. 
This equivalence, stated in the next theorem,  was first noted by Brimkov, Fast and Hicks in \cite{BFH19}.

\begin{thm} \label{t: forts}  {\rm \cite{BFH19}}
Let $G$ be a graph.  Then $B\subseteq V(G)$ is a zero forcing set if and only if $B$ intersects every fort.
\end{thm}

  A \emph{minimal fort} is a fort that is not properly contained in any other fort.   
  
 \begin{cor}\label{c: forts} Let $G$ be a graph.  Then $B\subseteq V(G)$ is a zero forcing set if and only if $B$ intersects every minimal fort.
 \end{cor}
 
\begin{rem}\label{r:minZminF} 
Let $G$ be a graph and let $B$ be a minimal zero forcing set.   Let $B'=\{u\in B: u\mbox{ is in some minimal fort of }G\}$. For every fort $F$,  $F$ contains a minimal fort $F'$.  Since $B$ is a zero forcing set, $ \emptyset\ne B\cap F'=B'\cap F'\subseteq B'\cap F$, so $B'$ is a zero forcing set.  Since $B$ is minimal, $B=B'$. So each element of $B$ is an element of some minimal fort of $G$.   
\end{rem}

\begin{rem}\label{r:fortunion}
In general, a superset of a fort $F$ need not be a fort, because if vertex $v$ is added to $F$ there may be a vertex $u$ such that $u$ has no neighbors in $F$ but is adjacent to $v$. However, there are conditions under which all or certain supersets of a fort are a fort.  If $F$ is a fort and  $N(v)\subseteq F$, then $F\cup\{v\}$ is a fort.
If $F$ is a fort such that  $N(v)\subseteq F$ for every $v\not\in F$, then every superset of $F$ is a fort.
A union of forts is a fort: Suppose $F_1$ and $F_2$ are forts of $G$.  If $u\not\in F_1\cup F_2$, then $u$ has zero or at least two neighbors in each of $F_i$, so $u$ has zero or at least two neighbors in $F_1\cup F_2$
\end{rem}

A private fort is a natural analog of a private neighbor and is used to define zero forcing irredundance. 

\begin{defn} 
Let  $G$ be a graph and let $S\subseteq V(G)$.  For $x\in S$, a fort $F$ of $G$ is a \emph{private fort} of $x$ (relative to $S$), if $S\cap F=\{x\}$.   A \emph{minimal private fort} of $x$ is a private fort of $x$ that does not properly contain another private fort of $x$.
\end{defn}

\subsection{Zero forcing irredundance}
Analogously to a zero forcing set that intersects every fort, we have that $D\subseteq V(G)$ is a dominating set if and only if $D$ intersects every closed neighborhood. With irredundance being a natural extension of domination, we introduce the concept of zero forcing irredundance. 

\begin{defn}\label{def: zir-set and number}
Let  $G$ be a graph and let $S\subseteq V(G)$.  Then $S$ is a \emph{$\Z$-irredundant set} or \emph{ZIr-set} if every element of $S$ has a  {private fort}.  
The \emph{lower ZIr number} is \[\zir(G)=\min\{|S|:S\mbox{ is a maximal ZIr-set}\}.\]
The \emph{upper ZIr number} is \[\ZIR(G)=\max\{|S|:S\mbox{ is a maximal ZIr-set}\}.\]
\end{defn}

We use the terms \emph{upper ZIR set} and \emph{lower zir set} to refer to maximal ZIr-sets $S$  and $S'$ such that $|S|=\ZIR(G)$ and $|S'|=\zir(G)$, respectively.

In light of Corollary \ref{c: forts}, it is reasonable to consider replacing private forts in Definition \ref{def: zir-set and number} with private \emph{minimal} forts. 
However, this modification does not align with the definition of irredundance: $v$ is a private neighbor of $u$ if $S\cap N[v]=\{u\}$ (where $N[v]=N(v)\cup\{v\}$), and this definition does not require that   $S\cap N[x]=\{u\}$ implies that $N[v]\subseteq N[x]$. It also has an undesired consequence: For every graph $G$ and every vertex $v\in V(G)$, the set $\{v\}$ is an irredundant set. The same statement for $\Z$-irredundance does not necessarily hold when restricting to private minimal forts. 

\subsection{Structure of the paper}
Basic results  about the relationship between $\Z, \zbar$, $\ZIR$, and $\zir$ are presented in Section \ref{s:prelim}, with many parallels to $\gamma, \Gamma$, IR, and ir.  Examples of upper ZIR and lower zir numbers for specific graph families are given  in Section \ref{s:families}.  Section \ref{s:bounds} presents bounds on upper ZIR and lower zir numbers and graphs having extreme upper ZIR and lower zir numbers are characterized in Section \ref{s:ext}. The effect of the join and corona operations is studied in Section \ref{s:graph-ops} and it is shown in Section \ref{s:noncompare} that  
the lower zir number is noncomparable to several lower bounds for the zero forcing number.
Section \ref{s:compute} provides information for future study of computational complexity and discusses the software used here for computations.
Finally, Section \ref{s:conclude} presents a summary table of upper ZIR and lower zir numbers, zero forcing numbers and other related parameters for various graph families and discusses directions for future research. 

\subsection{Notation}

We conclude this introduction with some additional basic graph theory notation and terminology that we will use. A graph $G=(V(G),E(G))$ is simple, undirected, finite,  and $V(G)\ne \emptyset$.
Standard symbols are used for well-known graph families: $K_n$ denotes a complete graph of order $n$; this graph has an edge between every pair of vertices. $\ol{K_n}$ denotes an empty graph of order $n$, which has no edges. $P_n$ and $C_n$ denote the path  graph and cycle graph of order $n$; the vertices of $P_n$ can be numbered $v_1,\dots,v_n$ so that the edges are $v_iv_{i+1}$ for $i=1,\dots,n-1$, and similarly for the cycle with the addition of edge $v_nv_1$. 
A disjoint union is denoted by $\du$ and is used for both sets and graphs (meaning the vertex sets are disjoint). The join $G\vee H$ of two disjoint graphs $G$ and $H$ has $V(G\vee H)=V(G)\cup V(H)$ and $E(G\vee H)=E(G)\cup E(H)\cup \hat E$ where $\hat E=\{uw:u\in V(G) \mbox{ and }w\in V(H)\}$. If $G$ and $H$ are disjoint and $G$ has order $n_G$,  the \emph{corona} of $G$ with $H$, denoted by $G\circ H$, is the graph obtained from the disjoint union of $G$ and $n_G$ copies of $H$ by joining the $i$th vertex of $G$ and vertices of the $i$th copy of $H$ for each $i=1,\dots,n_G$.  

Vertices $u$ and $w$ of a graph $G$ are \emph{twins} if $N(u)=N(w)$ (\emph{independent twins}) or $N[u]=N[w]$ (\emph{adjacent twins}). A set of vertices $\{w_1,\dots, w_r\}$ is a \emph{set of twins} if every two vertices in the set are independent twins or every two vertices in the set are adjacent twins. The next result is well known and useful in the study of zero forcing.

\begin{prop}\label{p:twins}{\rm \cite[Proposition 9.15]{HLS22}} 
If a graph has a set of twins $\{w_1,\dots, w_r\}$, then each zero forcing set must contain at least $r-1$ of these twins.
\end{prop}

 \section {Preliminary results} \label{s:prelim}

Many of the results in this section are not difficult to prove but provide intuition and establish a collection of tools we shall frequently use. Since zero forcing is one of the primary inspirations for this work, we begin by investigating its relationship with ZIr-sets.  We have an immediate parallel with dominating and irredundant sets.  

\begin{rem}\label{r:Z+ZIr}
Let $S\subseteq V(G)$   be both a zero forcing set and ZIr-set. Then we can see that $S$ is a minimal zero forcing set and maximal ZIr-set: Since every element $x\in S$ has a private fort $F_x$, $(S \setminus \{x\})\cap F_x=\emptyset$ and $S \setminus \{x\}$ is not a zero forcing set. Since a superset of a zero forcing set is a zero forcing set, and no proper superset of $S$ is miminal, $S$ is a maximal ZIr-set.
\end{rem}

\begin{prop} \label{c:zir-Z-Zbar-ZIR}
Let $G$ be a graph and $S\subseteq V(G)$.  Then $S$ is a minimal zero forcing set if and only if $S$ is a maximal  ZIr-set of $G$ and $S\cap F\ne\emptyset$ for every fort $F$ of $G$.
\end{prop}
\begin{proof}
Suppose that $S$ is a minimal zero forcing set. Let $x\in S$. Since $S$ is minimal, $S \setminus \{x\}$ is not a zero forcing set. By Theorem \ref{t: forts} there exists a fort $F_x$ such that $(S \setminus \{x\})\cap F_x=\emptyset$ and $S\cap F_x \ne\emptyset$. Thus $F_x$ is a private fort of $x$ relative to $S$. This argument holds for every vertex in $S$ and so $S$ is a ZIr-set. Theorem \ref{t: forts} guarantees $S\cap F\ne\emptyset$ for every fort $F$ of $G$ and hence $S$ is a maximal ZIr-set.

Now suppose that $S$ is a maximal ZIr-set of $G$ and $S\cap F\ne\emptyset$ for every fort $F$ of $G$. Then Theorem \ref{t: forts} implies that $S$ is a zero forcing set, so $S$ is a minimal zero forcing set by Remark \ref{r:Z+ZIr}. 
\end{proof}

Note that since $V(G)$ is a fort for any graph $G$, $1 \leq \zir(G)$.

\begin{cor}\label{cor: ZIr and ZF bounds}
Let $G$ be a graph. Then  $1\le \zir(G)\le \Z(G)\le\zbar(G)\le\ZIR(G)$.
\end{cor}

Not every maximal ZIr-set is a zero forcing set, as the next example shows.  

\begin{ex}
Consider the cycle $C_5$.   
The minimal forts are sets of three vertices that do not contain three consecutive vertices (consecutive in the order around the cycle).  
It follows that every  set of two vertices is a maximal ZIr-set, but not every set of two vertices is a zero forcing set since a zero forcing set of $C_5$ must contain two adjacent vertices.
\end{ex}

Next we present elementary results about ZIr sets.  
 As the next remark illustrates, it often suffices to restrict to connected graphs when studying ZIr-sets.

\begin{rem}\label{r:disconn-sum}
    If $G=G_1\du \dots \du G_k$, then $S$ is a ZIr-set of $G$ if and only if $S\cap G_i$ is a ZIr-set of $G_i$ for $i=1,\dots,k$.  Thus $\zir(G_1\du\dots\du G_k) = \zir(G_1) +\dots+ \zir(G_k)$ and $\ZIR(G_1\du\dots\du G_k) = \ZIR(G_1) +\dots+ \ZIR(G_k)$
    \end{rem}

 The next example is an immediate consequence.

\begin{ex}\label{e:empty}
Let $n\geq 1$. Since $\zir(K_1) = 1$, $\zir(\overline{K_n}) = Z(\overline{K_n}) = \overline{Z}(\overline{K_n}) = \ZIR(\overline{K_n}) = n$.
\end{ex}

\begin{rem}\label{r:extbds} Suppose $G$ is a graph of order $n\ge 2$ that has an edge.  Every fort in a connected component containing an edge must contain at least two vertices. Adding an additional private fort to a union of private forts must add a new vertex,  so the union of  all $k$ private forts associated with a maximal ZIr-set of $k$ vertices must contain at least $k+1$ vertices.  Thus $\ZIR(G)\le n-1$.
\end{rem}

\begin{rem}
Let $G$ be a graph with no isolated vertices.  Then every minimal zero forcing set is a maximal ZIr-set and no vertex is in every minimal zero forcing set \cite{smallparam}. Thus  there does not exist a vertex $v$ that is contained in every maximal ZIr-set of $G$.
\end{rem}

\begin{lem}\label{lem: complement of ZIR} 
Let $G$ be a  graph with no isolated vertices  
and let $S\subseteq V(G)$. If  $N[v]\subseteq S$ for some vertex $v$,  then $S$ is not a ZIr-set, or equivalently, if   $S$ is a ZIr-set of $G$, then $V(G)\setminus S$ is a dominating set of $G$.
\end{lem}
\begin{proof}
Suppose there exists a vertex $v\in S$ such that $N(v) \subseteq S$. Having assumed that $G$ has no isolated vertices, there exists a vertex $u\in N(v)$.  Let $F\subseteq V(G)$ such that $F\cap S = \{u\}$. Since $N(v)\subseteq S$, $v$ is adjacent to exactly one element of $F$. Thus there does not exist a private fort of $u$ relative to $S$, and so $S$ is not a ZIr-set of $G$.  \end{proof}

  It is well known (and easy to see) that $\delta(G)\le \Z(G)$.  The next  two results strengthen this bound.

 \begin{prop}\label{deg-d}
Let $G$ be a graph and $d\geq 1$. If $v_1,\ldots,v_d \in V(G)$ all have degree at least $d$, then $S = \{v_1,\ldots,v_d\}$ is a ZIr set.
\end{prop}
\begin{proof}
Assume $v_1,\ldots,v_d \in V(G)$ all have degree at least $d$. It suffices to show that  the set $F_i =(V(G)\setminus S)\cup \{v_i\}$ is a private fort of $v_i$ relative to $S$ for every $i = 1,\ldots,d$. This is certainly the case if $d=1$ since $V(G)$ is a fort, so suppose $d\geq 2$. Let $v\in V(G)\setminus F_i$ for some $i = 1,\ldots, d$. If $v$ has two or more neighbors in $V(G) \setminus S$, then $|N(v)\cap F_i| \geq 2$. So, assume $v$ has at most one neighbor in $V(G) \setminus S$. Since $v\in S$ and $|N(v)| \geq d$, $v$ is adjacent to exactly one vertex in $V(G) \setminus S$ and to $v_i$ because $S\setminus \{v\} \subseteq N(v)$. Thus $|F_i\cap N(v)|\geq 2$ as required.
\end{proof}

\begin{cor}\label{prop: min deg bnd}
  If  $G$ is a graph and $\delta(G)\ge 1$, then  every set of size $\delta(G)$ is a ZIr-set.  For every graph $G$, $\delta(G)\le \zir(G)$.
\end{cor}

The next example shows the bound in Corollary \ref{prop: min deg bnd}  is sharp.

\begin{ex}\label{p:clique} Let $n\geq {2}$. Then  $\zir(K_n)=n-1=\ZIR(K_n)$ by Corollary \ref{prop: min deg bnd} (since  $\delta(K_n)=n-1$) and Remark \ref{r:extbds}.  
\end{ex}
Example \ref{p:clique} implies  that  $\Z(K_n) =\overline{Z}(K_n) = n-1$, but this is already known \cite{AIM,Zrecon}.
In the next section we establish values of the parameters $\zir, \Z,\zbar,$ and $\ZIR$ for additional graph families.  These results show that both $\zir(G)=\Z(G)=\zbar(G)=\ZIR(G)$ and $\zir(G)<\Z(G)<\zbar(G)<\ZIR(G)$ are possible (Example  \ref{p:clique} and Proposition \ref{ex:all-diff}). 

 \begin{obs}\label{privatefortsum}
   Let $S$ be a  ZIr-set of a graph $G$,  of order $n$, let $x_i\in S$ for $i=1,\dots,k$, and let $F_i$ be a private fort of $x_i$.  Then $|S|\le n-|\cup_{i=1}^kF_i|+k$.
\end{obs}


\section{Determining ZIR and zir for graph families}\label{s:families}

In this section we establish the values of the upper and lower ZIr numbers for various families of graphs.  In doing so, we highlight various techniques that will be used throughout the paper. We also include known zero forcing numbers and upper zero forcing numbers in each of the following results as the connections between the parameters and the underlying graph is rather interesting. See Table \ref{table1} for a summary of the values of $\zir, \Z, \zbar$, and $\ZIR$ determined for various graphs in this and other sections.

Let $G$ be a graph on $n$ vertices. The next example illustrates that the gap between $\zir(G)$ and $\Z(G)$ (and hence also $\ZIR(G)$) can be as large as $n-3$. It is  known that $ \Z(K_{q,p})=q+p-2$ and $\zbar(K_{q,p}) = q+p-2$ \cite{AIM,Zrecon}.
 
\begin{ex}\label{ex:comp-pipart}\label{p:star}  Suppose $1\le q\le p$ and let $U = \{u_1,\ldots,u_q\}$ and $W = \{w_1,\dots,w_p\}$ be the partite sets of $K_{q,p}$.  Any set $S$ that omits at least one vertex $u_k$ from $U$ and at least one vertex $w_\ell$ from $W$ is a ZIr-set, since then  $\{u_i,u_k\}$ is a private fort of $u_i$ and $\{w_j,w_\ell\}$ is a private fort of $w_j$  (if $q=1$, the forts $\{u_i,u_k\}$ do not exist). If $S$ is a ZIr-set containing $U$, then $S=U$ because a private fort of $u_i$ necessarily contains $W$ and similarly for $W\subseteq S$  (note that when $q=1$, the private fort of $u_1$ is $V(K_{1,p})$). Thus $\Z(K_{q,p}) = \zbar(K_{q,p}) = \ZIR(K_{q,p})=q+p-2$ and 
$\zir(K_{q,p}) = q$.
\end{ex}

\begin{prop}\label{prop: cycle_zir_numbers}
Let $n \geq 4$. Then $\zir(C_n) = \Z(C_n) = \zbar(C_n) = 2$ and $\ZIR(C_n) = \lf\frac n2\rf$.
\end{prop}
\begin{proof}
It is known that $Z(C_n) = \overline{Z}(C_n) = 2$ \cite{AIM,Zrecon}.  Together with Corollary \ref{prop: min deg bnd}, this implies $\zir(C_n)=2$. Let $v_1,\ldots,v_n$ denote the vertices of $C_n$, where $v_n v_1 \in E(C_n)$ and $v_i v_{i+1}\in E(C_n)$ for $i=1,\ldots,n-1$.    
Let $S = \{v_{2i-1} : i = 1,\ldots, \lfloor \frac n2 \rfloor\}$ and $F_i = (V(C_n)\setminus S )\cup \{v_i\}$.  By construction, $F_{2i-1}$ is a private fort of $v_{2i-1}$ for $i=1,\dots,\lfloor \frac n2 \rfloor$. Thus $S$ is a ZIr-set. 
 Since $|S| = \lfloor \frac n2 \rfloor$, $\ZIR(C_n)\ge \lf\frac n 2 \rf$.

Since every subset of $V(C_n)$ that does not contain a pair of consecutive vertices has cardinality at most $\lfloor \frac n2\rfloor$,  it suffices to show that if $S$ is a ZIr-set of $C_n$ that contains a pair of adjacent vertices, then $ |S| = 2$.

Assume that $S$ is a ZIr-set of $C_n$ that contains a pair of adjacent vertices. Without loss of generality let $\{v_1,v_2\}\subseteq S$. Let $F$ be a fort of $C_n$. Suppose, to obtain a contradiction, that  $F\cap \{v_1,v_2\} = \emptyset$. Then $v_3 \notin F$ since  otherwise $|N(v_2)\cap F| = 1$. Recursively applying this argument implies $v_k\notin F$ for $k = 4,\ldots,n$ and hence $F = \emptyset$, a contradiction. Thus every fort of $C_n$ intersects $ \{v_1,v_2\}$ non-trivially. Since every element of $S$ must have a private fort, it follows that $ |S| = 2$.   
\end{proof}

Observe that $C_3=K_3$, so the parameters are determined by  Example \ref{p:clique} (and  $\ZIR(K_3) = 2 > \lfloor \frac32 \rfloor$).

\begin{prop}\label{p:path}
For $n\ge 1$, $\zir(P_n) = \Z(P_n) = 1$. For $n\ge 4$,  $\zbar(P_n)=2$. For $n\ge 5$, $\ZIR(P_n) = \lf\frac {n-1}2\rf$.
\end{prop}
\begin{proof}
It is well known that $Z(P_n) = 1$ for $n\ge 1$, which implies $\zir(P_n)=1$. For $n\ge 4$, $\overline{Z}(P_n) = 2$ \cite{Zrecon}. Assume $n\ge 5$.  Let $v_1,\dots,v_n$ denote the vertices of $P_n$, where $v_i v_{i+1}\in E(P_n)$ for $i=1,\dots,n-1$. Just as in the proof of Proposition \ref{prop: cycle_zir_numbers}, if $S$ is a ZIr-set of $P_n$ that contains a pair of adjacent vertices, then $|S| \le 2$. Observe that every fort in $P_n$ contains $v_1$ and $v_n$.

Let $S = \{v_{2i} : i = 1,\ldots, \lfloor \frac{n-1}2 \rfloor\}$ and $F_i = (V(P_n)\setminus S )\cup \{v_i\}$. Note that $|S| = \lfloor \frac{n-1}2 \rfloor$. By construction, $F_{2i}$ is a private fort of $v_{2i}$ for $i=1,\dots,\lfloor \frac{n-1}2 \rfloor$. Thus $S$ is a ZIr-set. Moreover, $S$ is maximal since $|S|\geq 2$ and the addition of any vertex to $S$ would introduce a pair of consecutive vertices or a vertex of degree 1. Since any set that does not contain a pair of consecutive vertices and does not contain $v_1$ or $v_n$ has cardinality at most $\lfloor \frac{n-1}2\rfloor$ we conclude that $\ZIR(P_n) = \lfloor \frac{n-1}2 \rfloor$.
\end{proof}

\begin{obs}\label{o:small-path}
Observe that  $P_2=K_2$ and $P_3=K_{1,2}$, so $\ZIR(P_2)=1$ and   $\ZIR(P_3)=1$ from Examples \ref{p:clique} and  \ref{p:star}.  It can be shown computationally \cite{sage} that  $\ZIR(P_4)=2$. 
\end{obs}

\begin{figure}[!h]
\centering
\scalebox{1}{
\begin{tikzpicture}[scale=2,every node/.style={draw,shape=circle,outer sep=2pt,inner sep=1pt,minimum size=.2cm}]		
\node[fill=none, label={[yshift=-27pt]$v_0$}]  (00) at (0,0) {};
\node[fill=none, label={[yshift=-5pt]$v_1$}]  (01) at (-0.342,0.94) {};
\node[fill=none, label={[yshift=-5pt]$v_2$}]  (02) at (0.342,0.94) {};
\node[fill=none, label={[yshift=-5pt]$v_3$}]  (03) at (0.985,-0.174) {};
\node[fill=none, label={[yshift=-15pt, xshift=12pt]$v_4$}]  (04) at (0.643,-0.766) {};
\node[fill=none, label={[yshift=-5pt]$v_6$}]  (05) at (-0.985,-0.174) {};
\node[fill=none, label={[yshift=-15pt, xshift=-12pt]$v_5$}]  (06) at (-0.643,-0.766) {};
		
\draw[thick] (00)--(01)--(02)--(00);
\draw[thick] (00)--(03)--(04)--(00);
\draw[thick] (00)--(05)--(06)--(00);
\end{tikzpicture}}
\caption{The friendship graph $\fr(3)$}
\label{fig:F3}\vspace{-7pt}
\end{figure}
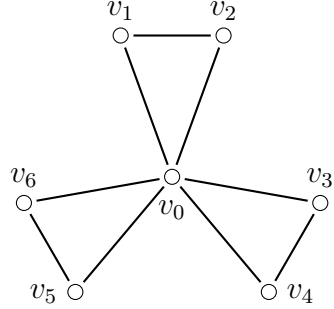

For $k\ge 2$, the $k$th-\emph{friendship graph} $\fr(k)$  is the graph with $2k+1$
vertices and $3k$ edges constructed as the join of one vertex to $k$ disjoint copies of $K_2$;   $\fr(3)$ is shown in Figure \ref{fig:F3}.  It is known that $\Z(\fr(k))=k+1$  \cite[Theorem 9.5]{HLS22}.\vspace{-5pt}

\begin{prop}\label{p:friend}
For $k\ge 2$, $\zir(\fr(k)) = \Z(\fr(k)) = \zbar(\fr(k)) = \ZIR(\fr(k)) = k+1$. 
\end{prop}

\begin{proof}
Let $k\geq 2$ and let $V(\fr(k)) = \{v_0,\ldots,v_{2k}\}$, where $v_0$ is the unique vertex of degree $2k$ and $v_j$ is adjacent to $v_{j+1}$ for $j\equiv 1 \mod 2$ and $1 \le j < 2k$.  We begin by examining the forts of $\fr(k)$.  For $i=1,\dots,k$, let $F_i=\{v_{2i-1},v_{2i}\}$ and observe that $F_i$ is a (minimal) fort. A union of forts of the form $F_i$ is a fort by Remark \ref{r:fortunion}. There are $2^k$ (minimal) forts that contain $v_0$ and exactly one element of each $F_i$ for $i=1,\dots,k$; a superset of such a fort is a fort by Remark \ref{r:fortunion}, and is called a \emph{center-fort}. These are the only forts of $\fr(k)$ because a fort that does not contain $v_0$ is a union of forts of the form $F_i$ and a fort that does contain $v_0$ is a center-fort.

There are only $k$ forts of the form $F_i$ so if $S$ is a ZIr-set and $|S|\ge k+1$, then one of the private forts relative to $S$ must be a center-fort $F'$, which contains at least $k+1$ vertices.  Thus $|V(\fr(k))\setminus F'|\le k$ and $|S|\le k+1$.   This implies $\ZIR(\fr(k))\le k+1$.

Let $S$ be a  ZIr-set. We show that $S$ is contained in a ZIr-set with cardinality $k+1$, which implies every maximal ZIr-set has cardinality $k+1$.  Suppose that  $S\cap F_{\ell}= \emptyset$ for some $\ell\in \{1,\dots,k\}$ and choose a set of private forts relative to $S$ for the vertices of $S$. Without loss of generality, we may assume that every chosen private fort that is not a center-fort is of the form $F_i$ for some $i\in \{1,\dots,k\}$.  We show that $S'=S\cup\{v_{2\ell-1}\}$ is also a ZIr-set with private fort $F_{\ell}$ for $v_{2\ell-1}$ and appropriate modifications of private forts of the vertices in $S$. Suppose  that the chosen private fort $F$ for some $v_j\in S$ contains $v_{2\ell-1}$. Note that $F$ is a center-fort. Since $v_j\not\in F_{\ell}$,    $F'=\lp F\cup \{v_{2\ell}\}\rp\setminus \{v_{2\ell-1}\}$ is a private fort relative to $S'$ for $v_j$. Thus $S$ is not maximal, and there is a ZIr-set $\wh S$ such that $S\subseteq \wh S$  and $\wh S\cap F_i\ne \emptyset$ for $i=1,\dots,k$. If $F_i\subseteq \wh S$ for some $i$, then $|\wh S|\ge k+1$. So assume $|\wh S\cap F_i|=1$ for $i=1,\dots,k$. If $v_0\in \wh S$, then $|\wh S|\ge k+1$. Suppose $v_0\not\in \wh S$. For $\wh S\cup \{v_0\}$, we have private forts $F_i$ for elements of $\wh S$ and the center-fort that does not contain any element of $\wh S$ for $v_0$.  Thus $\wh S\cup \{v_0\}$ is a ZIr-set and every maximal ZIr-set has cardinality at least $k+1$.
\end{proof}

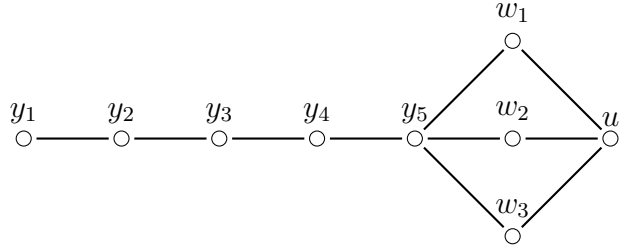
\begin{figure}[!h]
\centering
\scalebox{1}{
\begin{tikzpicture}[scale=1.3,every node/.style={draw,shape=circle,outer sep=2pt,inner sep=1pt,minimum size=.2cm}]		
\node[fill=none, label={[yshift=-5pt]$y_1$}]  (00) at (0,0) {};
\node[fill=none, label={[yshift=-5pt]$y_2$}]  (01) at (1,0) {};
\node[fill=none, label={[yshift=-5pt]$y_3$}]  (02) at (2,0) {};
\node[fill=none, label={[yshift=-5pt]$y_4$}]  (03) at (3,0) {};
\node[fill=none, label={[yshift=-5pt]$y_5$}]  (04) at (4,0) {};
\node[fill=none, label={[yshift=-5pt]$w_2$}]  (05) at (5,0) {};
\node[fill=none, label={[yshift=-5pt]$u$}]  (06) at (6,0) {};
\node[fill=none, label={[yshift=-5pt]$w_1$}]  (07) at (5,1) {};
\node[fill=none, label={[yshift=-5pt]$w_3$}]  (08) at (5,-1) {};
		
\draw[thick] (00)--(01)--(02)--(03)--(04)--(05)--(06)--(07)--(04);
\draw[thick] (04)--(08)--(06);
\end{tikzpicture}}
\caption{The graph $H(3,5)$, which has different values for $\zir,\Z,\zbar,$ and $\ZIR$}
\label{f:H35}
\end{figure}

For $r\ge 3$, denote the two vertices in the smaller partite set of $K_{2,r}$ by $u$ and $u'$ and the set of vertices in the larger partite set by $W=\{w_1,\dots,w_r\}$. For odd $s\ge 5$, denote the set of vertices of $P_s$ by $Y=\{y_1,\dots y_s\}$ in path order. Let $H(r,s)$ be the graph of order $r+s+1$ constructed from $K_{2,r}$ and $P_s$ by identifying  $u'$ with $y_s$, retaining the label $y_s$ (see Figure \ref{f:H35}). The next example shows that all the parameters $\zir,\Z,\zbar,$ and $\ZIR$ can be different (for example, $\zir(H(3,5))=2$, $\Z(H(3,5))=3$, $\zbar(H(3,5))=4$ and $\ZIR(H(3,5))=5$ \cite{sage}).

\begin{prop}\label{ex:all-diff}
  For $ r\ge 2$ and odd $s\ge 3$, $\zir(H(r,s))=2$, $\Z(H(r,s))=r$, and $\ZIR(H(r,s))= r+\frac{s-1}2$. For $ r\ge 2$ and odd $s\ge 5$,  $\zbar(H(r,s))=r+1$. For $ r\ge 2$, $\zbar(H(r,3))=r$.  
\end{prop}
\bpf Let $r\ge 2$ and let $s$ be odd with $s\ge 3$.
By  Proposition \ref{p:twins}, any zero forcing set must contain at least $r-1$ of the vertices $\{w_1,\dots, w_r\}$ since this is a set of twins. The only minimal zero forcing sets are $r-1$ of these vertices together with exactly one of the following sets: $\{u\}$, $\{y_1\}$, $\{y_s\}$, or  (if $s\ge 5$) $\{y_i,y_{i+1}\}$ for $i=2,\dots, s-2$. Thus  $\Z(H(r,s))=r$;   $\zbar(H(r,s))=r+1$ if $s\ge 5$ and $\zbar(H(r,3))=r$.

If $F$ is a fort such that $F\cap Y\ne \emptyset$, then $y_1\in F$ and at least every other vertex of $Y$ is in $F$. For the ZIr-set $S_1=\{u,y_1\}$, the only private fort $F_{u}$ of $u$ is $\{u,w_1,\dots,w_r\}$, so $S_1$ is maximal and $\zir(H(r,s))\le 2$. 
Since $S_1$ is a ZIr-set, $\{y_1\}$ is not maximal. Since every other vertex has degree at least two, every other vertex can also be in a ZIr-set of two vertices by Proposition \ref{deg-d}. Thus $\zir(H(r,s))= 2$.   

Let $D=\{y_1,y_3,\dots,y_{s-2},y_s,u\}$ and $S_2 = V(H(r,s))\setminus D$.  Every vertex $v\in S_2$ is adjacent to at least two vertices in $D$, so $\{v\}\cup D$ is a private fort of $v$ and $S_2$ is a ZIr-set. This shows that $\ZIR(H(r,s)\ge r+\frac{s-1}2$.  

Now let $S$ be a ZIr-set. There are $r+1$ vertices not in $Y$, and they cannot all be in $S$  because if $w_j\in S$ then a private fort of $w_j$ must contain $u$ or some $w_k$ with $k\ne j$. So $|S\cap \{u,w_1,\dots,w_r\}|\le r$. If $S$ does not contain two consecutive vertices in $Y$, then $|S\cap Y|\le \frac{s-1}2$. So assume that  $S$ contains two consecutive vertices $y_{i-1},y_i$. Then $S\cap Y\subseteq \{y_{i-1}, y_i \}$ because a fort $F$ such that $F\cap Y\ne\emptyset$ cannot omit two consecutive vertices of $Y$. Thus $|S\cap Y|\le 2\le \frac{s-1}2$. In all cases, $|S|\le  r+\frac{s-1}2$
\epf

The value of $\ZIR$ is established for additional families of graphs in Theorems  \ref{ZIR-n/2}, \ref{ZIR=2n/5} and \ref{ZIR-n/3} and Proposition \ref{p:wheel} below.


\section{Bounds}\label{s:bounds} 

In this section we use known results about 2-domination to establish bounds on $\ZIR(G)$ and other methods to establish additional bounds on $\ZIR(G)$. A \emph{$2$-dominating set} of a graph $G$ is a subset $D\subseteq V(G)$ such that every vertex not in $D$ is adjacent to at least two vertices in $D$. The minimum cardinality of a 2-dominating set is the \emph{$2$-domination number} of the graph, indicated by $\gamma_2(G)$. The idea used in the next proof to show that $ n-\gamma_2(G)\le \ZIR(G)$ was already used informally in Proposition \ref{ex:all-diff}  and other results.

\begin{prop}\label{p:2dom}
Let $G$ be a graph of order $n\ge 2$ with no isolated vertices. Then\break$ n-\gamma_2(G)\le \ZIR(G) \leq n-\gamma(G)$.  Both bounds are sharp.
\end{prop}
\bpf 
Let $S$ be $\ZIR$-set of $G$ such that $|S|=\ZIR(G)$. By Lemma~\ref{lem: complement of ZIR}, $V(G)\setminus S$ is a dominating set and therefore $\gamma(G)\leq |V(G)|-|S|=n-\ZIR(G)$.
   
Let $D$ be a $\gamma_2$-set of $G$ and let $S = V(G)\setminus D$. For $v\in S$ set $F_v = D\cup\{v\}$ is a private fort of $v$ relative to $S$. Therefore $S$ is a ZIr-set and $\ZIR(G)\geq |S|$.

Example \ref{ex:comp-pipart} shows $K_{q,p}$ with $q\ge 2$ realizes equality for the upper bound and Proposition \ref{prop: cycle_zir_numbers} shows $C_{n}$ realizes equality for the lower bound. 
\epf

Note that the  previous lower bound does not need the restriction that $G$ has no isolated vertices.

We use known results of Caro and Roddity \cite{caro90} and Cockayne et al.~\cite{cockayne85} bounding $\gamma_2(G)$ for graphs with specified minimum degree to bound $\ZIR(G)$ for such graphs; note that the results in \cite{caro90} and \cite{cockayne85} are more general and imply Theorems \ref{t:2dom:del3} and \ref{t:2dom:del2}. 

\begin{thm}\label{t:2dom:del3}{\rm \cite{caro90}}
Let $G$ be a graph of order $n$. If $\delta(G)\geq 3$, then $\gamma_2(G)\leq \frac{n}{2}$.
\end{thm}

\begin{cor}\label{c:ZIR2n2}
Let $G$ be a  graph of order $n\ge 2$ with $\delta(G)\geq 3$. Then $\ZIR(G) \geq \frac{n}{2}$.
\end{cor}

\begin{figure}[!h]
\centering
\scalebox{1}{
\begin{tikzpicture}[scale=1,every node/.style={draw,shape=circle,outer sep=2pt,inner sep=1pt,minimum size=.2cm}]

\node[fill=none, label={[yshift=-5pt]$a_{1}$}]  (1) at (-1.5,1.866) {};
\node[fill=none, label={[yshift=-25pt]$b_{1}$}]  (2) at (0,1) {};
\node[fill=none, label={[yshift=-5pt]$c_{1}$}]  (3) at (1.5,1.866) {};
\node[fill=none, label={[yshift=-5pt]$d_{1}$}]  (4) at (0,2.732) {};
\node[fill=none, label={[xshift=5pt, yshift=-5pt]$a_{2}$}]  (5) at (2.366,0.366) {};
\node[fill=none, label={[yshift=-5pt]$b_{2}$}]  (6) at (0.866,-0.5) {};
\node[fill=none, label={[yshift=-25pt]$c_{2}$}]  (7) at (0.866,-2.232) {};
\node[fill=none, label={[xshift=10pt, yshift=-7.5pt]$d_{2}$}]  (8) at (2.366,-1.366) {};
\node[fill=none, label={[yshift=-25pt]$a_{3}$}]  (9) at (-0.866,-2.232) {};
\node[fill=none, label={[yshift=-5pt]$b_{3}$}]  (10) at (-0.866,-0.5) {};
\node[fill=none, label={[xshift=-5pt, yshift=-5pt]$c_{3}$}]  (11) at (-2.366,0.366) {};
\node[fill=none, label={[xshift=-10pt, yshift=-7.5pt]$d_{3}$}]  (12) at (-2.366,-1.366) {};
		
\draw[thick] (2)--(3)--(4)--(2)--(1)--(4);
\draw[thick] (6)--(7)--(8)--(6)--(5)--(8);
\draw[thick] (10)--(11)--(12)--(10)--(9)--(12);
\draw[thick] (1)--(11);
\draw[thick] (3)--(5);
\draw[thick] (7)--(9);
\end{tikzpicture}}
\caption{The graph $N_3$}
\label{f:N3}\vspace{-7pt}
\end{figure}
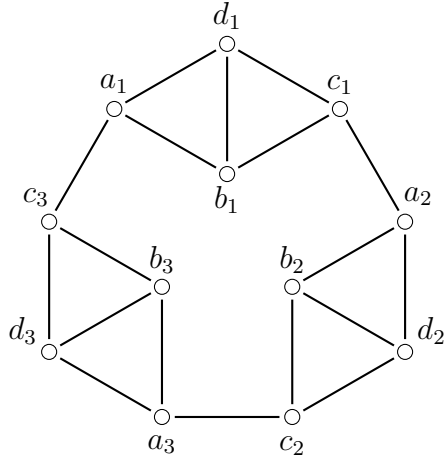

This bound is sharp, as is shown in the next result. A \emph{diamond} is a $K_4$ with one edge deleted. Let $D_1,D_2,\dots, D_n$ be $n$ diamonds each with vertex set $V_i=\{a_i,b_i,c_i,d_i\}$ where the edge $a_ic_i$ is missing. Define the graph consisting of the $k$ diamonds together with the edges $c_ia_{i+1}$ for $1\leq i\leq k-1$ and $c_ka_1$  to be the $k$th \emph{necklace}, denoted by $N_k$. The 3rd necklace $N_3$ is illustrated in Figure \ref{f:N3}.\vspace{-5pt}

\begin{thm}\label{ZIR-n/2}
Let $k\ge 2$. Then $\ZIR(N_k)=2k = \frac12 |V(N_k)|$.
\end{thm}
\bpf
By Corollary~\ref{c:ZIR2n2}, $\ZIR(N_k)\geq 2k$.
Let $S$ be a ZIr-set of $N_{k}$. If $|V_i\cap S|\leq 2$ for all $i$, then $|S|\leq 2k$.  
Note that $|S\cap V_i|=4$ is impossible because $V_i=N[b_i]\not \subseteq S$ by Lemma \ref{lem: complement of ZIR}. We show that if $|S\cap V_i| = 3$, then $|S\cap V_{i-1}|,|S\cap V_{i+1}|\leq 1$.  Once this is established, we conclude that $|S|\le 2k$.

First, we assume that $c_{j-1}$ and $a_j$ cannot be in a private fort of the vertices in $V_j\cap S$, and show $|V_j\cap S|\leq 1$. Let $x\in \{b_j,c_j,d_j\}\cap S$ and let $F_x$ be a private fort of $x$. If $F_x\cap \{b_j,d_j\}\neq \varnothing$, then $a_j$ must be adjacent to two vertices in the fort and $c_{j-1}\not\in F_x$, so $\{b_j,d_j\}\subseteq F_x$. If $x=c_j$, then $b_j\in F_x$ or $b_j$ needs two neighbors in $F_x$ so $d_j\in F_x$, so $\{b_j,c_j,d_j\}\subseteq F_x$.  
Therefore $|V_j\cap S|\leq 1$. By symmetry, if $c_{j-1}$ and $a_j$ cannot be in a private fort of the vertices in $V_{j-1}\cap S$, then $|S\cap V_{j-1}|\leq 1$. Thus
\begin{linenomath}
\begin{align}\label{eq:Nk} 
c_{j-1},a_j \mbox{ cannot be in a private fort of any }v\in V_j\cap S&\implies  |V_j\cap S|\leq 1;\\
\label{eq:Nk2} c_{j-1},a_j\mbox{ cannot be in a private fort of any }v\in V_{j-1}\cap S&\implies |V_{j-1}\cap S|\leq 1.
\end{align}
\end{linenomath}
Now suppose that $|V_i\cap S|= 3$. First assume that $a_i,b_i,c_i\in S$. Note that $d_i$ cannot belong to a private fort of a vertex not in $V_i$, since then $b_i$ would be adjacent to only one vertex in the fort. Because $d_i$ (and $b_i$) cannot be in a private fort of a vertex not in $V_i$, $c_{i-1}$ cannot belong to any private fort of a vertex in outside $V_i$ (or $a_i$ would have only one neighbor in this fort). Thus $a_{i}$ and $c_{i-1}$ cannot belong to any private fort of a vertex in $V_{i-1}$,  and so $ |S\cap V_{i-1}|\leq 1$. By symmetry, $|S\cap V_{i+1}|\leq 1$. 

Now assume that $b_i, c_i,d_i\in S$. Note that $a_i$ cannot belong to a private fort of any vertex outside $V_i$,  
since then $b_i$ and $d_i$ would be adjacent to only one vertex in that fort. Furthermore, if $c_{i-1}$ is in a private fort of a vertex not in $V_i$, $a_i$ would only be adjacent to one vertex in that fort. Similarly, $a_{i+1}$ cannot belong to a private fort of a vertex not in $V_i$, since then $c_i$ would be adjacent to only one vertex in that fort. By \eqref{eq:Nk} and \eqref{eq:Nk2}, $|S\cap V_{i-1}|\leq 1$ and $|S\cap V_{i+1}|\leq 1$. 
\epf\vspace{-7pt}

\begin{rem}\label{r:neck} 
It is known that $\Z(N_k)=k+2$  \cite[Theorem 9.5]{HLS22}.  Suppose $S$ is a minimal zero forcing set. Since $b_i$ and $d_i$ are  twins, at least one of $b_i$ and $d_i$ must be in $S$  by  Proposition \ref{p:twins}. In order for the forcing to start, there must be a blue vertex with three blue neighbors so there must be at least two other vertices in $S$, and with these vertices we have a zero forcing set.  Thus $\zbar(N_k)=k+2$.  
\end{rem} 

\begin{thm}\label{t:2dom:del2}{\rm  \cite{cockayne85}} 
Let $G$ be a graph of order $n$. If $\delta(G) = 2$, then $\gamma_2(G)\leq \frac{2n}{3}$.
\end{thm}

\begin{cor}\label{c:ZIRn3}
Let $G$ be a graph of order $n\ge 3$  and $\delta(G) = 2$. Then $\ZIR(G) > \frac{n}{3}$.
\end{cor}
\bpf
From Proposition~\ref{p:2dom} and Theorem~\ref{t:2dom:del2} it follows that $\ZIR(G) \geq \frac{n}{3}$. Graphs with 2-domination number equal to $ \frac {2n} 3$ were characterized in \cite{favaron08} as $G = H\circ K_2$ for any graph $H$. We show that $\ZIR(G) >\frac{n}{3}$.

Label the vertices of $H$ as $v_1,v_2,\dots, v_{n_H}$ and the vertices of $K_2$ joined to vertex $v_i$ as $a_i$ and $b_i$. The set $S=\{v_1\}\cup \{a_i : i=1,\dots,n_H\}$ is a ZIr-set of $G$ where the private fort of each $a_i$ is $\{a_i,b_i\}$ and the private fort of $v_1$ is $ V(G)\setminus \{a_i : i=1,\dots,n_H\}$.
\epf

Theorem \ref{ZIR-n/3} provides a family of graphs $G$ satisfying $\ZIR(G) = \frac13 |V(G)|$,  but these graphs have $\delta(G) = 1$,  so we do not know  how close to the bound $\frac13 |V(G)|$  we can get for graphs $G$ with $\delta(G) = 2$. 

\begin{figure}[!h]
\centering
\scalebox{1}{
\begin{tikzpicture}[scale=1,every node/.style={draw,shape=circle,outer sep=2pt,inner sep=1pt,minimum size=.2cm}]		
\node[fill=none, label={[label distance=-5pt]30:$v_{1,3}$}]  (00) at (1.2377,1.5979) {};
\node[fill=none, label={[label distance=-5pt]$v_{1,2}$}]  (01) at (0,0.6986) {};
\node[fill=none, label={[label distance=-5pt]120:$v_{1,1}$}]  (02) at (-1.2377,1.5979) {};
\node[fill=none, label={[yshift=-7.5pt]$v_{1,5}$}]  (03) at (-0.7649,3.0528) {};
\node[fill=none, label={[yshift=-7.5pt]$v_{1,4}$}]  (04) at (0.7649,3.0528) {};
\node[fill=none, label={[xshift=10pt, yshift=-10pt]$v_{2,3}$}]  (05) at (0.7649,-1.8708) {};
\node[fill=none, label={[yshift=-7.5pt]$v_{2,2}$}]  (06) at (0.605,-0.3493) {};
\node[fill=none, label={[xshift=10pt, yshift=-10pt]$v_{2,1}$}]  (07) at (2.0026,0.2729) {};
\node[fill=none, label={[xshift=10pt, yshift=-10pt]$v_{2,5}$}]  (08) at (3.0263,-0.864) {};
\node[fill=none, label={[xshift=-5pt, yshift=-7.5pt]$v_{2,4}$}]  (09) at (2.2614,-2.1889) {};
\node[fill=none, label={[xshift=-10pt, yshift=-10pt]$v_{3,3}$}]  (10) at (-2.0026,0.2729) {};
\node[fill=none, label={[yshift=-7.5pt]$v_{3,2}$}]  (11) at (-0.605,-0.3493) {};
\node[fill=none, label={[xshift=-10pt, yshift=-10pt]$v_{3,1}$}]  (12) at (-0.7649,-1.8708) {};
\node[fill=none, label={[xshift=5pt, yshift=-7.5pt]$v_{3,5}$}]  (13) at (-2.2614,-2.1889) {};
\node[fill=none, label={[xshift=-10pt, yshift=-10pt]$v_{3,4}$}]  (14) at (-3.0263,-0.864) {};
		
\draw[thick] (00)--(01)--(02)--(03)--(04)--(00)--(07)--(08)--(09)--(05)--(06)--(07);
\draw[thick] (05)--(12)--(13)--(14)--(10)--(11)--(12);
\draw[thick] (02)--(10);
\end{tikzpicture}}
\caption{The graph $H_3$}
\label{fig:ZIR=n-gamma2}\vspace{-7pt}
\end{figure}
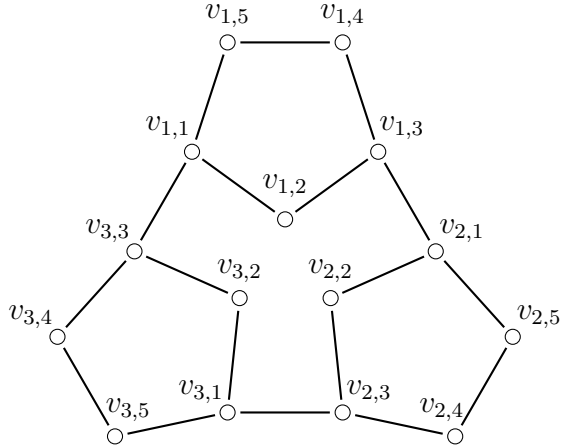

The next family of graphs realizes $\ZIR(G) = n-\gamma_2(G)$ and $\delta(G)=2$. Let $H_k$ be the graph consisting of $k$ 5-cycles, where the vertices of the $i$-th cycle is labeled $v_{i,1},v_{i,2},\dots, v_{i,5}$, together with the edges $v_{i,3}v_{i+1,1}$ for $i\in \{1,\dots,k-1\}$ and the edge $v_{k,3}v_{1,1}$; see Figure \ref{fig:ZIR=n-gamma2}. \vspace{-5pt}

\begin{thm}\label{ZIR=2n/5}
The graph $H_k$ of order $n=5k$ has $\delta (H_k)=2$,  $\gamma_2(H_k)=\frac{3n}{5}$, and $ \ZIR(H_k)=\frac{2n}{5}=n-\gamma_2(H_k)$.   Furthermore, $\Z(H_k)=k+2$.
\end{thm}

\bpf
Let $B_i$ denote the set of vertices of the $i$-th cycle of $H_k$. First we show that $\gamma_2(H_k)=\frac{3n}5$. Let $D$ be a 2-dominating set of $H_k$. If $v_{i,2}\in D$, then to 2-dominate $v_{i,4}$ and $v_{i,5}$ $|D\cap B_i\setminus\{v_{i,2}\}|\geq 2$. It therefore follows that $|D\cap B_i|\geq 3$. On the other hand, if $v_{i,2}\not\in D$, then certainly $v_{i,1}, v_{i,3}\in D$. Since $v_{i,4}$ and $v_{i,5}$ are not 2-dominated yet, it follows that $|D\cap\{v_{i,4},v_{i,5}\}|\geq 1$ and so $|D\cap B_i|\geq 3$. Hence $\gamma_2(H_k)\geq  \frac{3n}5$. Since $\{v_{i,2}, v_{i,4}, v_{i,5} : i= 1,\dots,k\}$ is a 2-dominating set, equality follows.

Since $\ZIR(H_k)\ge \frac{2n}5$, it suffices to show that $\ZIR(H_k)\le \frac{2n}5$. Let $S$ be a ZIr-set of $H_k$. If $|S\cap B_i|\leq 2$ for every $i\in \{1,\dots,k\}$ then the result follows. Assume therefore that $|S\cap B_i|\geq 3$ for some $i\in \{1,\dots,k\}$. If $|S\cap B_i|>3$, there exists a  degree-2 vertex  $v$ in $B_i$ such that $N[v]\subseteq S$, which is impossible according to Lemma~\ref{lem: complement of ZIR}. Therefore $|S\cap B_i|=3$.  To complete the proof, we show that $|S\cap B_{i-1}|\leq 1$ and $|S\cap B_{i+1}|\leq 1$.  Let $F_{i,j}$  denote a private fort of $v_{i,j}$ (relative to $S$).

First we  examine $B_i\cap S$ further,  showing $v_{i,1}$ and $v_{i,3}$ cannot both be in $S$.  To obtain a contradiction, assume $\{v_{i,1},v_{i,3}\}\subseteq S$. Since $N(v_{i,2})=\{v_{i,1},v_{i,3}\}$, $v_{i,2}\in F_{i,1}$. If $v_{i,5}\in F_{i,1}$, then $v_{i,4}$ would necessarily also be in $F_{i,1}$, contradicting $|S\cap B_i|=3$. So $v_{i,5}\not\in F_{i,1}$, and it follows that $v_{i,4}\in F_{i,1}$ to ensure that $v_{i,5}$ is adjacent to at least two vertices in $F_{i,1}$. A similar argument for a private fort $F_{i,3}$ shows that $v_{i,5}\in F_{i,3}$. In this case, $B_i \subseteq F_{i,1}\cup F_{i,3}$. Thus $|B_i\cap S|=2$, which is a contradiction. Hence $|\{v_{i,1},v_{i,3}\}\cap S|\leq 1$.
 
For the two cases $|\{v_{i,1},v_{i,3}\}\cap S|=  0$ and $|\{v_{i,1},v_{i,3}\}\cap S|= 1$, we prove both the following statements: 
\ben[$(a)$]
\item
There are $j$ and $j'$ such that  $v_{i+1,1}\in F_{i,j}$, and  $v_{i-1,3}\in F_{i,j'}$.
\item
$v_{i+1,1}\not\in F_{i+1,j}$ for any $j$ and $v_{i-1,3}\not\in F_{i-1,j}$ for any $j$. 
\een

\noi\underline{Case $|\{v_{i,1},v_{i,3}\}\cap S|= 0$}: 
 
First, we  prove  
(a) for both $v_{i+1,1}$ and $v_{i-1,3}$. Since $|B_i\cap S|=3$, $B_i\cap S=\{v_{i,2}, v_{i,4},v_{i,5}\}$. Then $v_{i,3}$ is the only neighbor of $v_{i,4}$ not in $S$ and therefore $v_{i,3}\not\in F_{i,2}$. To ensure that $v_{i,3}$ is adjacent to two vertices of $F_{i,2}$, $v_{i+1,1}\in F_{i,2}$. Similarly, $v_{i-1,3}\in F_{i,2}$.

To prove (b), suppose $v_{i+1,1}\in F_{i+1,j}$ for some $j\in \{1,\dots,5\}$. Since both of the neighbors of $v_{i,3}$ in $B_i$ are in $S$, $v_{i,3}$ is necessarily in $F_{i+1,j}$ but then $v_{i,4}$ is not adjacent to two vertices of $F_{i+1,j}$. Therefore $v_{i+1,1}\not\in F_{i+1,j}$ for any $j$ and similarly $v_{i-1,3}\not\in F_{i-1,j}$ for any $j$.
\smallskip

\noi\underline{Case $|\{v_{i,1},v_{i,3}\}\cap S|=1$}:

Assume that $\{v_{i,1},v_{i,3}\}\cap S=\{v_{i,3}\}$. It follows from Lemma~\ref{lem: complement of ZIR} that $N[v_{i,4}]\not\subseteq S$ and  thus at most one of $v_{i,4}$ or $v_{i,5}$ is in $S$. Since $|B_i\cap S|=3$, $v_{i,2}\in S$  and exactly one of $v_{i,4}$ or $v_{i,5}$ is in $S$.
 
First, we prove (a) for $v_{i+1,1}$. Either $v_{i,4}\in S$ or $v_{i,4}\in F_{i,5}$. In both cases $v_{i,3}$ is adjacent to a vertex in $F_{i,j}$ where $j=4$ or $5$, and therefore $v_{i+1,1}\in F_{i,j}$. 
 
Next, we prove (a) for $v_{i-1,3}$.   Let $j'=4$ or $5$, according to which of $v_{i,4}$ or $v_{i,5}$ is in $S$. If $v_{i,1}\in F_{i,j'}$, then  $v_{i,2}$ would have only one neighbor in $F_{i,j'}$.   Thus $v_{i,1}\not \in F_{i,j'}$, and  $v_{i-1,3}\in F_{i,j'}$ to ensure that $v_{i,1}$ is adjacent to two vertices of $F_{i,j'}$. 
  
To prove (b), suppose $v_{i+1,1}\in F_{i+1,j}$ for some $j\in \{1,\dots,5\}$.  Since $v_{i,2},v_{i,3}\in S$, $v_{i,4}$ is necessarily in $F_{i+1,j}$ and therefore $v_{i,5}\in S$. To ensure that $v_{i,5}$ is adjacent to two vertices of $F_{i+1,j}$, $v_{i,1}\in F_{i+1,j}$ but then $v_{i,2}$ is only adjacent to one vertex of $F_{i+1,j}$. Therefore $v_{i+1,1}\not\in F_{i+1,j}$ for any $j$. Now suppose that $v_{i-1,3}\in F_{i-1,j}$ for some $j\in \{1,\dots,5\}$. Since $v_{i,2},v_{i,3}\in S$, $v_{i,1}\not\in F_{i-1,j}$ and therefore $v_{i,5}\in F_{i-1,j}$. But then $v_{i,4}\in S$ and only adjacent to one vertex of $F_{i-1,j}$. Hence $v_{i-1,3}\not\in F_{i-1,j}$ for any $j$.

\smallskip

To complete the proof we show that $|S\cap B_{i+1}|\leq 1$ and $|S\cap B_{i-1}|\leq 1$.

From (a), $v_{i+1,1}\in F_{i,\ell}$ for some specific $\ell\in \{1,\dots,5\}$ and it therefore follows that $|\{v_{i+1,2},v_{i+1,3}\}\cap F_{i,\ell}|\geq 1$ and $|\{v_{i+1,4},v_{i+1,5}\}\cap F_{i,\ell}|\geq 1$. 

Now consider the vertices in $S\cap B_{i+1}$. From (b) it follows that $v_{i+1,1}\not \in F_{i+1,j}$ for any $j\in \{1,\dots,5\}$, which implies $v_{i+1,1}\not \in S$. Also we can see that $v_{i,3}\not\in F_{i+1,j}$ either:  If $v_{i,3}\in S$ then clearly $v_{i,3}\not\in F_{i+1,j}$.  So assume $v_{i,3}\not \in S$ and $v_{i,3}\in F_{i+1,j}$.   If $v_{i,1}\not\in S$, then $v_{i,2}, v_{i,4},v_{i,5}\in S$ and $v_{i,4}$ is adjacent to only one vertex in $F_{i+1,j}$. If $v_{i,1}\in S$, then also $v_{i,2}\in S$ and $v_{i,2}$ is adjacent to only one vertex in $F_{i+1,j}$.  Thus $v_{i,3}\not\in F_{i+1,j}$.

Since $v_{i,3}\not\in F_{i+1,j}$ for any $j$, $\{v_{i+1,2},v_{i+1,5}\}\not\subseteq S$,  otherwise $v_{i+1,1}$ is adjacent to only one vertex in $F_{i+1,2}$ and $F_{i+1,5}$.

Suppose first that $\{v_{i+1,2},v_{i+1,5}\}\cap S=\{v_{i+1,2}\}$. Then $v_{i+1,5}\in F_{i+1,2}$ to ensure that $v_{i+1,1}$ is adjacent to two vertices of $F_{i+1,2}$. Recall $v_{i+1,1}\in F_{i,\ell}$, so $v_{i+1,3}\in F_{i,\ell}$ to ensure that $v_{i+1,2}$ is adjacent to two vertices of $F_{i,\ell}$. If also $v_{i+1,4}\in S$, then $v_{i+1,5}\in F_{i+1,4}$. This is impossible since $v_{i+1,1}$ is then adjacent to only one vertex of $F_{i+1,4}$.  Thus $\{v_{i+1,2},v_{i+1,5}\}\cap S=\{v_{i+1,2}\}$ implies $|S\cap B_{i+1}|\leq 1$. 
If $\{v_{i+1,2},v_{i+1,5}\}\cap S=\{v_{i+1,5}\}$ it follows similarly that $|S\cap B_{i+1}|\leq 1$. 

Now, if $\{v_{i+1,2},v_{i+1,5}\}\cap S=\emptyset$, then $|S\cap B_{i+1}|\leq 2$. Assume $|S\cap B_{i+1}|=2$ with $v_{i+1,3}, v_{i+1,4}\in S$. Then $v_{i+1,2},v_{i+1,5}\in F_{i, \ell}$  since $v_{i+1,1}\in F_{i, \ell}$ and $v_{i+1,3}, v_{i+1,4}\not \in F_{i,\ell}$. But then $v_{i+1,4}$ is adjacent to only one vertex in $F_{i,\ell}$, a contradiction. It follows that $|S\cap B_{i+1}|\leq 1$. It follows similarly that $|S\cap B_{i-1}|\leq 1$ and hence $\ZIR(H_k)\geq \frac{2}{5}n$ as required.

For the zero forcing number, observe that $\{v_{1,3},v_{1.4}\}\cup\{v_{i,2}:i=1,\dots,k\}$ is a zero forcing set of size  $k+2$. Since $\{v_{i,2},v_{i,4},v_{i,5}\}$ is a fort,  a zero forcing set must contain at least one vertex of each cycle. Since $\deg(v_{i,1}) = \deg(v_{i,3})=3$, forcing cannot move from one cycle to the next unless some cycle has three vertices or two adjacent cycles have two each.  Thus $\Z(H_k)=k+2$. 
\epf

Finally, we present bounds that do not involve $\gamma_2(G)$.

\begin{rem}
Note that $\Z(G)\le\ZIR(G)$, and it is known that $\Z(G(n,p))=n-o(n)$ where $G(n,p)$ is the Erd\H{o}s-R\'enyi random graph of order $n$ with edge probability $p$ \cite[Section 9.11]{HLS22}.
\end{rem}

\begin{prop}\label{prop: max deg bnd}
Let $G$ be a connected graph of order $n\ge 2$.  Then $\ZIR(G)\le\frac{\Delta(G)}{\Delta(G)+1}n$  and this bound is sharp.
\end{prop}
\bpf
Let $S$ be a maximal ZIr-set of $G$. Let $X$ denote the set of edges incident to a vertex in $S$ and a vertex in $V\setminus S$.   Since $G$ is connected, a vertex $x\in S$ has a neighbor $w$.  If $w\in S$, then $x$ has another neighbor $u\in F_{w}$ (the private fort of $w$), and $u\not \in S$.  Thus every vertex in $S$ is adjacent to at least one vertex not in $S$ and therefore $|X|\geq |S|$.
Furthermore, every vertex in $V\setminus S$ is adjacent to at most $\Delta(G)$ vertices of $S$ and therefore $\Delta(G)|V\setminus S|\geq |X|$. It follows that 
\[
\Delta(G)(n-|S|)= \Delta(G)|V\setminus S|\geq |X|\geq |S|.
\]
Since $S$ is any maximal ZIr-set, $\ZIR(G)\le\frac{\Delta(G)}{\Delta(G)+1}n$.  By Example \ref{p:clique}, $\ZIR(K_n)=n-1=\frac{n-1}{(n-1)+1}n$. 
\epf

\begin{cor}\label{c:Delta-bds}
If $G$ is a connected cubic graph of order $n\ge 4$, then  $\frac 1 2 n \le \ZIR(G)\le \frac 3 4 n$ and the lower bound is sharp.
\end{cor}
\bpf
The inequalities are immediate from    Corollary \ref{c:ZIR2n2} and Proposition \ref{prop: max deg bnd}.  The lower bound is sharp by Theorem \ref{ZIR-n/2}.
\epf

The only example we know where the upper bound is equality is $\ZIR(K_4)=3$. 


\section{Extreme values of $\zir$ and $\ZIR$}\label{s:ext}

 In this section we characterize graphs  of order $n$ having the extreme high values $n$ and $n-1$ for $\zir$ and $\ZIR$ and characterize graphs having the extreme high value $n-2$  for $\zir$.  We also characterize graphs having the extreme low value 1 for $\zir$ and present a family of graphs with the lowest known $\ZIR$ value for arbitrarily large $n$.

\begin{rem}
Let $G$ be a graph of order $n$. It is well known that  $\Z(G)=n\Leftrightarrow G\cong \ol{K_{n}}\Leftrightarrow \zbar(G)=n$, and this  extends to $\zir$ and $\ZIR$: Since $\zir(G)=n$ implies $\Z(G)=n$ and $\ZIR(G)=n$ implies $G\cong \ol{K_n}$ by Remark \ref{r:extbds}, we have $G  \cong \ol{K_n}\Leftrightarrow \zir(G) = n \Leftrightarrow \ZIR(G)=n$.
\end{rem}

For a graph on $n$ vertices, it is well known that  $\Z(G)=n-1$ if and only if $G\cong K_{n-r}\du rK_1$ with $n-r\ge 2$. It was shown recently by Brimkov and Carlson in \cite{BC22} that $\zbar(G) = n-1$ if and only if $G\cong K_{n-r}\du rK_1$. We extend this to $\zir$ and $\ZIR$.

\begin{prop}\label{prop: ZIR le n-1}
Let $G$ be a graph of order $n\ge 2$.  
The following are equivalent:
\ben
\item\label{item:Kn} $G\cong K_{n-r}\du rK_1$ with $n-r\ge 2$. 
\item\label{item:zir=n-1} $\zir(G)=n-1$.
\item\label{item:Z=n-1} $\Z(G)=n-1$.
\item\label{item:Zbar=n-1} $\zbar(G)=n-1$.
\item\label{item:ZIR=n-1} $\ZIR(G)= n-1$.
\een
\end{prop}
\bpf 
If $G$ does not have an edge, then $G\cong \ol{K_n}$ and $\zir(\ol{K_n}) = \Z(\ol{K_n}) = \zbar(\ol{K_n}) = \ZIR(\ol{K_n}) = n$. 
So assume $G$ has an edge, so that $\ZIR(G)\le n-1$ by Remark \ref{r:extbds}.  Then with Example  \ref{p:clique} and  Corollary \ref{cor: ZIr and ZF bounds}, we have that \eqref{item:Kn} $\Rightarrow$ \eqref{item:zir=n-1} $\Rightarrow$ \eqref{item:Z=n-1} $\Rightarrow$ \eqref{item:Zbar=n-1} $\Rightarrow$ \eqref{item:ZIR=n-1}. Since $\zbar(H)=|V(H)|-1$ implies $H\cong K_{n-r}\du rK_1$ with $n-r\ge 2$ for any graph $H$ of order $n$ \cite{BC22},   \eqref{item:Zbar=n-1} $\Rightarrow$ \eqref{item:Kn}. Thus \eqref{item:Kn} -- \eqref{item:Zbar=n-1} are equivalent.

So assume $\ZIR(G)= n-1$. Then there is a ZIr-set $S$ such that $|S|=n-1$.  Since every set of $n-1$ vertices is a zero forcing set of $G$,  
$S$ is a minimal zero forcing set  by Remark \ref{r:Z+ZIr}.  Thus $\zbar(G)=n-1$.
\epf

The next result provides many examples of graphs with $\ZIR(G)=|V(G)|-2$.

\begin{cor}\label{c:HveeK2}
If  $G=H\vee 2K_1$ or $G=H\vee {K_2}$ with $H\not\cong K_{|V(H)|}$, then $\ZIR(G)=|V(H)|=|V(G)|-2$.  If $G=H\vee {K_2}$ and $H$ has no isolated vertices, then $\Z(G)=\Z(H)+2$.  
\end{cor}
\bpf 
The statement about $\ZIR$ is immediate from the previous proposition and Proposition \ref{p:2dom} (because $V(G)\setminus V(H)$ is a 2-dominating set). It is known  that $\Z(H\vee K_1)=\Z(H)+1$ for a graph $H$ with no isolated vertices \cite[Proposition 9.16]{HLS22}. Since $H\vee K_2\cong (H\vee K_1)\vee K_1$, this implies  $\Z(H\vee K_2)=\Z(H)+2$ if $H$ has no isolated vertices.  \epf

The previous corollary does not list $\Z(H\vee 2K_1)$ because $\Z(H)+1\le \Z(H\vee 2K_1)\le Z(H)+2$ and both bounds are sharp  (consider $H \cong K_r$ for the lower bound and $H\cong P_r$ for the upper bound). It is not true that all graphs $G$ of order $n$ with $\ZIR(G)=n-2$ have $\gamma_2(G)=2$, as seen in the next example.

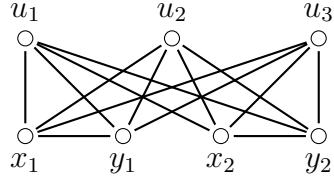
\begin{figure}[!h]
\centering
\scalebox{1}{
\begin{tikzpicture}[scale=1.3,every node/.style={draw,shape=circle,outer sep=2pt,inner sep=1pt,minimum size=.2cm}]		
\node[fill=none, label={[xshift=0pt, yshift=-5pt]$u_1$}]  (1) at (0,1) {};
\node[fill=none, label={[xshift=0pt, yshift=-5pt]$u_2$}]  (2) at (1.5,1) {};
\node[fill=none, label={[xshift=0pt, yshift=-5pt]$u_3$}]  (3) at (3,1) {};
\node[fill=none, label={[xshift=0pt, yshift=-25pt]$x_1$}]  (4) at (0,0) {};
\node[fill=none, label={[xshift=0pt, yshift=-25pt]$y_1$}]  (5) at (1,0) {};
\node[fill=none, label={[xshift=0pt, yshift=-25pt]$x_2$}]  (6) at (2,0) {};
\node[fill=none, label={[xshift=0pt, yshift=-25pt]$y_2$}]  (7) at (3,0) {};
		
\draw[thick] (1)--(4)--(5)--(2)--(6)--(7)--(3)--(6)--(1)--(5);
\draw[thick] (1)--(7)--(2)--(4)--(3)--(5);
\end{tikzpicture}}
\caption{A graph $G$ with $\ZIR(G)=|V(G)|-2$ and $\gamma_2(G)=3$.}
\label{fig:P3duP4comp}
\end{figure}

\begin{ex}\label{3K1vee2K2} 
Denote the vertices of the partite sets of $K_{3,4}$ by $U=\{u_1,u_2,u_3\}$ and $\{x_1,x_2,y_1,y_2\}$. Construct $G$ from $K_{3,4}$ by adding the edges $x_1y_1$ and $x_2y_2$. See Figure \ref{fig:P3duP4comp}. Note that no set of two vertices 2-dominates $G$ since the degree of every vertex is four. Since $U$ is a 2-dominating set, $\gamma_2(G)= 3$. We see that $S = \{u_1,u_2,u_3,y_1,y_2\}$ is a ZIr-set by noting that $F_{u_i} = \{u_i,x_1,x_2\}$ is a private fort of $u_i$ for $i=1,2,3$ and $F_{y_j}=\{y_j,x_j\}$ is a private fort of $y_j$ for $j=1,2$. Thus $\ZIR(G) = 5 = |V(G)|-2$. 
\end{ex}

There is a known characterization of $\Z(G)\ge n-2$ for graphs $G$ of order $n\ge 3$. In  \cite{AIM} it is shown that $\Z(G) \ge n - 2$ if and only if  
$G$ does not contain $P_4$, $P_3 \du K_{2}$, ${\rm dart}$, $\ltimes$, or $3K_2$ as an induced subgraph. In \cite{BHL04} it is shown that $G$ does not contain $P_4$, $P_3 \du K_{2}$, ${\rm dart}$, $\ltimes$, or $3K_2$ as an induced subgraph if and only if \[G \cong \overline{\lp K_{s_1}\du\dots \du K_{s_t}\du K_{q_1,p_1}\du K_{q_2,p_2}\du\cdots\du K_{q_k,p_k}\rp\vee K_r}\] with $s_i\ge 1$, $p_i\ge 1$, $q_i\ge 0$,  $t,k,r\ge 0$, and $t+k+r\ge 1.$ 
Observe that $K_1\cong K_{0,1}$ and $K_2\cong K_{1,1}$, so we may assume $s_i\ge 3$.  Without loss of generality we may also assume $p_i\ge q_i$ (which implies $p_i\ge 1$) and $q_i\ge q_{i+1}$. Furthermore, $K_{0,a_1}\du \dots \du K_{0,a_\ell}\cong K_{0,a_1+\dots+a_\ell}$, so we may assume $q_{k-1}\ge 1$. These results can be combined as in the next theorem. 

\begin{thm}\label{t:Zn-2} {\rm\cite{AIM, BHL04}} For a graph $G$ of order $n\ge 3$, $\Z(G)\ge n-2$ if and only if 
\[
\overline G\cong \lp K_{s_1}\du\dots K_{s_t}\du K_{q_1,p_1}\du K_{q_2,p_2}\du\cdots\du K_{q_k,p_k}\rp\vee K_r
\] 
with $s_i\ge 3$, $p_i\ge q_i$, $q_i\ge q_{i+1}$, $p_i\ge 1$, $q_{k-1}\ge 1$ and $q_k\ge 0$, $t,k,r\ge 0$, and $t+k+r\ge 1$. 
\end{thm}

 Note that $\Z(G)\ge |V(G)|-2$ implies $\ZIR(G)\ge |V(G)|-2$.  For example, the graph $G$ in Example \ref{3K1vee2K2} is covered by Theorem \ref{t:Zn-2} since $\Gc=K_3\du K_{2,2}$.  However, $\ZIR(G)\ge |V(G)|-2$ does not imply $\Z(G)\ge |V(G)|-2$, as seen in the next example.

\begin{ex}\label{rK2vee2K1}
Consider the graph $G = rK_2 \vee 2K_1$ with $r\ge 3$. Corollary \ref{c:HveeK2} establishes $\ZIR(G) = 2r = |V(G)|-2$. To see that $\Z(G)\le r+2$, note that a set consisting of one vertex from each $K_2$ and the vertices $u$ and $u'$ of  $2K_1$ is a zero forcing set of cardinality $r+2$, so $\Z(G)\le r+2 < 2r = |V(G)|-2$. Furthermore, $\Z(G)= r+2$, because any zero forcing set must contain at least one vertex of each $K_2$ and at least one of $u,u'$, and one more vertex (either contain both $u$ and $u'$ or contain both vertices of some $K_2$).
\end{ex}

We use 
 Theorem \ref{t:Zn-2} to characterize graphs $G$ having $\zir(G)=|V(G)|-2$.  

\begin{thm} 
Let $G$ be a graph of order $n\ge 3$. Then $\zir(G)= n-2$ if and only if 
\[
\overline G\cong {\lp  K_{q_1,p_1}\du K_{q_2,p_2}\du\cdots\du K_{q_k,p_k}\rp\vee K_r}
\] 
with $p_i\ge q_i$, $q_i\ge q_{i+1}$, $p_i\ge 1$, $q_{k-1}\ge 1$ and $q_k\ge 0$,  $k\ge 1$, $r\ge 0$, and  
\[
q_1\ge 2 \ \mbox{ or } \ ( q_1= 1 \mbox{ and } k\ge 2).
\]
\end{thm}
\begin{proof} 
Note that $\zir(G)=n-2$ implies $\Z(G)\ge n-2$.  Thus it suffices to assume $G$ has the form in Theorem \ref {t:Zn-2}, and determine which of these forms have $\zir(G)=n-2$. Recall $\zir(G) = n$ if and only if $G\cong nK_1$, and $\zir(G) = n-1$ if and only if $G\cong K_{n-r}\du rK_1$ with $n-r\ge 2$. This eliminates several cases: $t=1, k=0$ means $G\cong (s_1+r)K_1$; $t=0, k=1, q_1=0$ means $G\cong K_{p_1}\du rK_1$;  $t=0, k=1, q_1=1$ means $G\cong K_{p_1}\du (r+1)K_1$.  Thus $t=1, k=0$ and $t=0, k=1, q_1=0,1$ imply $\zir(G)>n-2$.  
This means 
$G$ has one of the following forms: \vspace{-5pt}  
\ben[(a)] 
\item\label{case:a} $t\ge 2$.\vspace{-5pt}
\item\label{case:b} $t= 1$ and $k\ge 1$. \vspace{-5pt}
\item\label{case:d} $t=0$  and $q_1\ge 2$.\vspace{-5pt}
\item\label{case:c} $t=0$, $k\ge 2$, and $q_1=1$.\vspace{-5pt}
\een

Observe that the presence or absence of $r$ isolated vertices in $G$ ($\vee K_r$ in $\ol G$) has no effect on whether or not $\zir(G)=n-2$, because each isolated vertex is its own private fort. Thus we assume $r=0$.

We first assume we have one of cases \eqref{case:a} and \eqref{case:b} and show that 
$\zir(G)<n-2$. In both these cases, $t\ge 1$ and $s_1\ge 3$. Let $U$ denote the vertices of $\ol{K_{s_1}}$. Since $t\ge 2$ or $k\ge 1$, $G$ is  connected. Let  
$S=V(G)\setminus U$; note
$S\ne \emptyset$. Since $G[U]$ is the empty graph on $s_1$ vertices and every vertex in $U$ is adjacent to every vertex not in $U$, $U$ is in every private fort of a vertex in $S$. Therefore $S$ is a maximal ZIr-set, where $\{v\}\cup U$ is a private fort of $v\in S$  with $|S|=n-|U|<n-2$.

Case \eqref{case:d} and \eqref{case:c} 
satisfy the description of $G$ in the statement of the theorem, so it remains to show that every graph $G$ of this form satisfies $\zir(G)= n-2$.   Denote the partite sets as $Z_1,\dots,Z_{2k}$, or $Z_{2k-1}$ if  $q_k=0$, where $|Z_{2j-1}|=p_j$ and $|Z_{2j}|=q_j$ for $j=1,\dots,k$. 
Let $S$ be a  maximal ZIr-set. We show  that $|S|\ge n-2$, which implies $\zir(G)=n-2$. 
The only case in which a partite set is a dominating set is if $q_k=0$; in this case each vertex of $Z_{2k-1}$ dominates $G$. Suppose first that $Z_{2k-1}$ is a dominating set and $\lp V(G)\setminus S\rp \subseteq Z_{2k-1}$.  If $| V(G)\setminus S|\le 2$, then $|S|\ge n-2$.  If $|V(G)\setminus S|\ge 3$ then $S$ would not be maximal, because for any distinct $u,w\in Z_{2k-1}$, $V\setminus \{u,w\}$ is a ZIr-set since $\{v,u,w\}$ is a private fort of $v$. So assume    $\lp V(G)\setminus S\rp \not \subseteq Z_{2k-1}$ or $Z_{2k-1}$ is not a dominating set. For all $i$, if $Z_i$ is not a dominating set of $G$, then $\lp V(G)\setminus S\rp \not\subseteq Z_i$  by Lemma \ref{lem: complement of ZIR}. This means $\lp V(G)\setminus S\rp \not\subseteq Z_i$ for all $i$, so there exist distinct indices $j$ and $\ell$ and vertices $z_j\in \lp V(G)\setminus  S\rp \cap Z_j$ and $z_\ell\in \lp V(G)\setminus  S\rp \cap Z_\ell$.  We show that $V\setminus \{z_j,z_\ell\}$ is a ZIr-set, so the maximality of $S$ implies $V\setminus \{z_j,z_\ell\}\subseteq S$ and $|S|\ge n-2$. Let $v\ne z_j, z_\ell$ and $v\in Z_i$.  If $i\ne j,\ell$, the $\{v,z_j,z_\ell\}$ is a private fort of $v$. If $i=j$, then $\{v,z_j\}$ is a private fort of $v$, and the case $i=\ell$ is similar.
\end{proof}

\begin{rem}\label{r:singleton}
Suppose $\zir(G)=1$ and $\{x\}$ is a maximal ZIr-set.  If $F_x$ is a fort that contains $x$, then $F\subseteq F_x$ for every fort $F$ that does not contain $x$, because if $F\not\subseteq F_x$, then there exists $y\in F\setminus\ F_x$ and $\{x,y\}$ is a ZIr-set.
\end{rem}

\begin{thm}\label{t:zir1}
Let $G$ be a graph of order $ n\ge 1$.  Then $\zir(G)=1$ if and only if $G = P_n$ or $G = K_{1,n-1}$.
\end{thm}
\bpf  
If $G=P_n$ or $G=K_{1,n-1}$, then $\zir(G)=1$ by Proposition \ref{p:path} or Example  \ref{p:star}. Assume $\zir(G)=1$. If $G$ had more than one connected component, then $\zir(G)\ge 2$ by Remark \ref{r:extbds}, so $G$ is connected. If $n=1$ then $G=P_1$.  So assume $n\ge 2$ and  $G\ne K_{1,p}$ with $p\ge 1$. Then $G$ has at least two vertices of degree at least two.  

Since $\zir(G)=1$, there exists a maximal ZIr-set $\{x\}$.  By Proposition \ref{deg-d}, 
$\deg x\not\ge 2$ since $\{x\}$ is maximal and there are two vertices of degree two or more.  Thus $\deg x=1$.

Now suppose that $G$ has a vertex of degree at least three. Let $c$ be a vertex such that $\deg c\ge 3$ and $\deg w\ge 3$ implies $\dist(x,w)\ge \dist(x,c)$. Label the vertices on the unique path from $c$ to $x$ by $ v_1,\dots,v_k$ with $\dist(v_i,c) = i$.  Since $c$ is a cut-vertex and $\deg c\ge 3$,  $F_1 = V(G)\setminus\{c,v_1,\dots,v_k,x\}$ is a fort.  Since $x\not\in F_1$ and $\{x\}$ is a maximal ZIr-set, it follows from Remark~\ref{r:singleton} that every fort that contains $x$ must also contain $F_1$. For any vertex $u$, $\deg u\ge 2$ implies $F_2 = V(G)\setminus\{u\}$ is a fort that contains $x$.  If $u\ne v_i$, then $F_1\not\subseteq F_2$.  
Thus $\deg u\ge 2$ implies $u = v_i$ for some $i$. Since $\deg c\ge 3$, there are two vertices $w_1, w_2\in N(c)$ with $w_1,w_2\ne v_1$.  Thus $\deg w_1 = \deg w_2=1$.  Thus $F_3 = V(G)\setminus\{c,w_1\}$ is a fort that contains $x$. But $F_1 \not\subseteq F_3$, which is a contradiction. Thus $G$ cannot have a vertex of degree greater than two. Since $\deg x=1$, $G$ is a path.
\epf

Note that $\ZIR(G)=1$ implies $\Z(G)=1$ and thus $G \cong P_n$. But $\ZIR(P_n)=1$ only when $n\leq 3$  (see Proposition~\ref{p:path} and Observation \ref{o:small-path}).

We do not know what the lowest value of $\ZIR(G)$ is over all graphs of order $n$ as $n\to\infty$.  A computer search of graphs of order $n\le 9$ found several graphs $G$ with order 9 and $\ZIR(G)=3$, including $C_3\circ 2K_1$. This led to Theorem \ref{ZIR-n/3} below  that $\ZIR(C_r\circ 2K_1)=\frac 1 3|V(C_r\circ 2K_1)|$, the lowest known value of $\ZIR(G)$ among graphs $G$ of order $n$ as $n$ becomes arbitrarily large. Before stating Theorem \ref{ZIR-n/3}, we introduce some notation and prove a lemma. We adopt the following labeling for vertices of $ H\circ tK_1$ where $r$ is the order of $H$: Label the vertices of $H$ as $u_0,\ldots,u_{r-1}$ and let the set of leaves adjacent to $u_i$ be denoted by $L_i=\{v_{i,j}:j= 1,\ldots, t\}$ for $i = 0,\ldots, r-1$. Furthermore, let $W_i=\{u_i\}\cup L_i$ for $i = 0,\ldots, r-1$.

\begin{lem}\label{l:circ-leaf} 
Let $H$ be a connected graph of order $r\ge 3$, let $ t\ge 2$, let $ G=H\circ tK_1$, and let $S$ be a  ZIr-set of $G$. Then $G$ has a ZIr-set $S'$ such that $|S|=|S'|$ and every vertex of $S'$ is a leaf.  
\end{lem}
\bpf
Let $S$ be a ZIr-set and let $S_i=S\cap W_i$.  Note that if $F$ is a fort and  $u_i\in F$, then $W_i\subseteq F$.  So  $u_i\in S$ implies  $v_{i,j}\not\in S$ for $j=1,\dots, t$ because none of these vertices could have a private fort relative to $S$.  
If $u_i\in S$, we wish to replace $u_i$ by an element of $L_i$ in $S$ to obtain a ZIr-set of the same cardinality with one fewer vertex of $H$. Repeating this process as needed will establish the result. To replace $u_i$ by  $v_{i,1}$, we need to verify that we can adjust other private forts if needed so that $v_{i,1}$ is not in the private fort of any element of $S$ except $u_i$.  Let $x\in S_k$ with $k\ne i$.  Choose a minimal private fort $F_x$. 
We show that minimality implies $|F_x\cap L_i|\le 1$: Note that $u_i\not\in F_x$ since $u_i\in S$. The only reason for any $v_{i,j}$ to be in $ F_x$ is that $u_i$ has exactly one neighbor in $F_x$ that is not in $L_i$. Thus including one leaf in $L_i$ is sufficient. If $F_x\cap L_i\ne \emptyset$, without loss of generality we modify $F_x$ so that $F_x\cap L_i =\{v_{i,t}\}$. Define $S' = S\setminus \{u_i\}\cup \{v_{i,1}\}$. Then $S'$ is a ZIr-set with private forts $\{v_{i,1}, v_{i,t}\}$ for $v_{i,1}$ and $F_x$ for $x\not\in W_i$. 
\epf

Note that maximality need not be preserved in moving to an all-leaf ZIr-set of the same cardinality.

\begin{thm}\label{ZIR-n/3}
Let $r\ge 3$. Then \[\Z(C_r\circ 2K_1) = \zbar(C_r\circ 2K_1) = \ZIR(C_r\circ 2K_1) = r = \frac 1 3 |V(K_r\circ 2K_1)|.\]
\end{thm}
\bpf
In $G = C_r\circ 2K_1$, we modify the labeling convention to label the two leaf neighbors of $u_i$ by $x_i$ and $y_i$. Index arithmetic is done in $\ZZ_{r}$, i.e., modulo $r$.  Observe that $x_i$ and $y_i$ are twins for $i = 1,\dots,r$, so at least one of them must be a zero forcing set. Thus $r\le \Z(G)\le  \zbar(G) \le\ZIR(G)$. 

Now we show that $\ZIR(G)\le r$. Let $W_i = \{u_i,x_i,y_i\}$. Let $S$ be a ZIr-set and let $S_i = S\cap W_i$. By Lemma \ref{l:circ-leaf}, we may replace any occurrence of $u_i$ in $S$ by $x_i$ and assume $S$ consists entirely of leaves. Thus the possible values of $s = |S_i|$ are $0$, $1$, or $2$. For $s = 0,1,2$ define $T_s = \{i: |S_i| = s\}$ and observe that $r = |T_2|+|T_1|+|T_0|$ and $|S| = 2|T_2|+|T_1|$.  We show that $|T_2|\le |T_0|$, which implies that $|S|\le r$. 
For each $z\in S$ we choose a minimal private fort $F_{z}$. If $z\in S_i$, $u_j\in F_{z}$, and $j\ne i$, then $j\in T_0$ because $u_j\in F_{z}$ implies $W_j\subseteq F_{z}$.

If $T_2=\emptyset$, then there is nothing to prove, so without loss of generality (by renumbering if necessary) assume  $0\in  T_2$, so $S_0=\{x_0,y_0\}$. Note that $y_0\not\in F_{x_0}$, which implies $u_0\notin F_{x_0}$. Thus $u_0$ needs at least two neighbors in $F_{x_0})$, so one of $u_{r-1}$ and $u_{1}$ must be in $F_{x_0}$. Without loss of generality, assume that $u_{1}\in F_{x_0}$, which implies $1\in T_0$. We proceed in increasing order around the cycle. If there are two consecutive cycle vertices $u_{j},u_{j+1}\not \in F_{x_0}$, then we stop as soon as this is encountered and set $j_0 = j$. If we do not stop because there is no $j$ such that $u_{j},u_{j+1}\not \in F_{x_0}$, then at least half of the cycle vertices $u_j$ are in $F_{x_0}$.  Thus $|T_0|\ge \frac r 2\ge |T_2|$. 

For the remainder of the proof we assume   $F_{z}$ omits two consecutive cycle vertices for each $z\in S$. The assumption that $u_{1}\in F_{x_0}$ and we have stopped with $u_{j_0},u_{j_0+1}\not \in F_{x_0}$  and $u_{j_0-1}\in F_{x_0}$ implies at least one of $x_{j_0}$ or $y_{j_0}$ is in $F_{x_0}$. 
Define $I_0=\{0,1,\dots,j_0-1\}$. Observe that $|I_0 \cap T_0|\ge  |I_0\cap T_2|$ and $j_0\not\in T_2$.

We start with the interval $I = I_0$ and expand as needed through the process described next, such that the interval $I = [0,b]$ retains the following properties (which are true for $I_0$):
\ben
\item\label{c1} $|I \cap T_0|\ge  |I\cap T_2|$.
\item\label{c3} If $b\in T_2$, then $b-1\in I\cap T_0$.  If $b\in T_0$, then $b+1\not\in  T_2$.
\een
If $T_2\subset I$, then the proof is complete by condition \eqref{c1}. If not, let $p$ be the first element of $T_2$ that is not in $I$ starting at $b$ and proceeding in increasing order. We examine a private fort $F_{x_p}$.  As in the construction of $I_0$, one of $u_{p-1}$ or $u_{p+1}$ is in $F_{x_p}$. If $u_{p-1}\in F_{x_p}$, then update the interval to $I = [0,\dots,p]$ and observe that the conditions \eqref{c1} and \eqref{c3} remain true because $p-1\in T_0$ and $j\not\in T_2$ for $j = b+1,\dots,p-2$ by the way $p$ was chosen. So suppose $u_{p-1}\not\in F_{x_p}$, which implies $u_{p+1}\in F_{x_p}$. We proceed in increasing order around the cycle. If for all $j=p+2,\dots,r-2$, $u_{j}\not \in F_{x_p}$ implies $u_{j+1} \in F_{x_p}$, then $I=\{0,\dots,r-1\}$ satisfies \eqref{c1} and the proof is complete. So let $j$ be the first index (in increasing order) such that $u_{j},u_{j+1}\not \in F_{x_p}$. Then, update the interval to $I = [0,\dots,j-1]$ and observe that the conditions \eqref{c1}--\eqref{c3} remain true. This completes the proof.
\epf

\section{Effect of graph operations on upper ZIr number} \label{s:graph-ops}

In this section we examine the effect of graph operations in ZIr-sets and ZIr numbers, primarily the upper ZIr number. We first consider the effect of removing a cut-vertex.

\begin{prop}\label{p:cutvtx} 
Let $G$ be a  connected graph such that $c$ is a cut-vertex of $G$ and $G-c$ has $\ell\ge 3$ components, denoted  by $H_1,\dots, H_\ell$ and ordered so that $\ZIR(H_1)\ge \dots\ge \ZIR(H_\ell)$.  For $i = 1,\dots,\ell$, let $S_i$ be a ZIr-set for $H_i$.  Then $\bigcup_{i\ne k} S_i$ is a ZIr-set for $k = 1,\dots,\ell$ and $\ZIR(G)\ge \ZIR(H_1)+ \dots+\ZIR(H_{\ell-1})$.   This bound is sharp.
\end{prop}
\bpf    
Fix $k\in\{1,\dots,\ell\}$ and let $S = \bigcup_{i\ne k} S_i$. For each $x\in S$, there is an $i\ne k$ such that $x\in S_i$. Then $x$ has a private fort  $F$ in $H_i$. If $|F\cap N(c)| = 0$ or $|F\cap N(c)|\ge 2$, then $F$ is a private fort of $x$ in G relative to $S$.  If  $|F\cap N(c)| = 1$, then $F' = F\cup V(H_k)$ is a private fort of $x$ in G relative to $S$. Since every $x\in S $ has a private fort of $x$ in G relative to $S$, $S$ is a ZIr-set, and by choosing $k = \ell$ we have $\ZIR(G)\ge \ZIR(H_1)+ \dots+\ZIR(H_{\ell-1})$. 

To see that the bound is sharp, consider the graph $G = C_r\circ 2K_1$.  Each cycle vertex is a cut-vertex with components $H_1  
 \cong  P_{r-1}\circ 2K_1$, $H_2  \cong K_1$, and $H_3  \cong K_1$.   By Theorem \ref{ZIR-n/3}, $\ZIR(C_r\circ 2K_1) = r$.  Furthermore, $\Z(P_{r-1}\circ 2K_1) = r-1$ since $P_{r-1}\circ 2K_1$ is a tree. Thus  $\ZIR(C_r\circ 2K_1) = r\ge \ZIR(H_1)+\ZIR(H_2) = \ZIR(P_{r-1}\circ 2K_1) + \ZIR(K_1)\ge \Z(P_{r-1}\circ 2K_1) + \ZIR(K_1) = (r-1)+1 = r $.
\epf

We now consider the join of two graphs $G\vee H$. Proposition \ref{prop: Join of G and H} handles the more general case and Proposition \ref{p:ZIR-GveeK1} is needed for the special case of $H = K_1$.

\begin{prop}\label{prop: Join of G and H}
Let $G$ and $H$ be graphs on $n_G\geq 2$ and $n_H\geq 2$ vertices, respectively. Then $n_G + n_H - 4 \leq \ZIR(G\vee H) \leq n_G + n_H - 1 $. Both bounds are sharp.  In particular, $\ZIR(P_{n_G}\vee P_{n_H})=n_G+n_H-4$.
\end{prop}
\begin{proof}
Let $D\subseteq V(G \vee H)$ such that $|D\cap V(G)| = 2$ and $|D\cap V(H)| = 2$. Then $D$ is a 2-dominating set. By Proposition~\ref{p:2dom}, $\ZIR(G \vee H) \geq n_G + n_H - 4$. Since $G\vee {H}$ has an edge, Remark~\ref{r:extbds} implies $\ZIR(G \vee H)\leq n_G + n_H - 1$.  

The upper bound is realized by $G \cong K_{n_G}$ and $H \cong K_{n_H}$. We now show that the lower bound is realized by $G  \cong P_{n_G}$ and $H  \cong P_{n_H}$ for $n_G\geq 7$ and $n_H\geq 7$. Let $v_1,\ldots, v_{n_G}$ and $u_1,\ldots,u_{n_H}$ denote the vertices of $G$ and $H$, respectively, where $v_i$ (respectively $u_i$) is adjacent to $v_j$ (respectively $u_j$) if and only if $i = j \pm 1$. Assume, by way of contradiction, that $\ZIR(G\vee H) \geq n_G + n_H - 3$. Then $G\vee H$ has a ZIr-set $S$ of size $n_G + n_H - 3$. 

First, suppose $V(G)\subseteq S$. Observe that every private fort of every vertex in $S\cap V(H)$ relative to $S$ is a fort of $H$, and therefore contains $u_1$ and $u_{n_H}$. Thus, $u_1,u_{n_H}\notin S$ and so there exists some $j\in\{2,n_G-1\}$ such that $u_j\in S$. This is a contradiction since every fort in $H$ containing $u_j$ has size at least 5, which  implies  
$|S|\le n_G + n_H - 5 + 1$. Thus $V(G)$ is not a subset of $S$. Similarly, $V(H)$ is not a subset of $S$.

Without loss of generality, assume $|V(G)\setminus S| = 1$ and $|V(H)\setminus S| = 2$. Since $n_H\geq 7$, there exists some $u_j\in V(H)$ such that $N[u_j]\cap V(H) \subseteq S$, and there exists some $u_k\in S\cap V(H)$ such that $u_k$ is not adjacent to $u_j$ and $j \not= k$. Let $F_{u_k}$ be a private fort of $u_k$ relative to $S$. Since every fort of $H$ is a dominating set of $H$, $F_{u_k}$ is not a fort of $H$. Thus, $F_{u_k}\cap V(G) = \{v_\ell\}$  where $V(G)\setminus S=\{v_\ell\}$. 
This is a contradiction since $N[u_j] \cap F_{u_k} = \{v_\ell\}$. Thus, $\ZIR(G\vee H) = n_G + n_H - 4$.
\end{proof}

The lower bound in the previous result can fail when $n_G = 1$ and $n_H = 2$, as illustrated by the friendship graph $\fr(k) = K_1\vee k  K_2$ since $|V(\fr(k))| = 2k+1$ and $\ZIR(\fr(k)) = k+1$ by Proposition \ref{p:friend}.

\begin{prop}\label{p:ZIR-GveeK1}
Let $G$ be a graph of order $n_G$  with no isolated vertices. Then 
\beq\label{eq:bds-gamma}
n_G - \gamma(G)\leq \ZIR(G\vee K_1) \leq n_G-\gamma(G)+1.
\eeq 
If $\ZIR(G) = n_G - \gamma(G)$, then $\ZIR(G\vee K_1) = n_G-\gamma(G)+1 = \ZIR(G)+1$. Both bounds  in \eqref{eq:bds-gamma} are sharp.
\end{prop}
\begin{proof} 
Let $u$ denote the vertex of $K_1$, let $S$ be a maximal ZIr-set of $G$, and let $D$ be a dominating set of $G$. Then $D\cup\{u\}$ is a 2-dominating set of $G\vee K_1$, so $\gamma_2(G\vee K_1)\le \gamma(G) +1$ and $\ZIR(G\vee K_1)\geq  (n_G+1) - (\gamma(G)+1)= n_G-\gamma(G)$ by Proposition~\ref{p:2dom}.

 Let $\hat S$ be a $\ZIR$-set of $G\vee K_1$. If $u$ is in a private fort $F$ of some vertex in $\hat S$, then $F$ contains at least $\gamma(G)$ vertices of $G$ since $u$ is adjacent to every vertex in $G$. Therefore 
$|\hat S|\leq (|V(G)|+1)-|F|+1 = n_G+1-(\gamma(G)+1)+1 = n_G-\gamma(G)+1.$
If $u$ is not in any private fort of a vertex in $\hat S$, then $u\not\in \hat S$ and every private fort relative to $\hat S$ in $G\vee K_1$ is also a private fort relative to $\hat S$ in $G$. Therefore $\hat S$ is a ZIr-set of $G$, and $\ZIR(G\vee K_1)=|\hat S|\leq \ZIR(G)\leq n_G-\gamma(G)$ by Proposition~\ref{p:2dom}. Thus in either case $\ZIR(G\vee K_1)\le n_G-\gamma(G)+1$.

By Lemma \ref{lem: complement of ZIR}, $V(G)\setminus S$ is a dominating set of $ G$ (recall $S$ is a maximal ZIr-set of $G$). Thus, $\{u\} \cup (V(G)\setminus S)$ is a private fort of $u$ in $G\vee K_1$ relative to $S\cup \{u\}$. So $\ZIR(G\vee K_1) \geq \ZIR(G) + 1$.  Thus $\ZIR(G) =n_G - \gamma(G)$ implies 
$\ZIR(G\vee K_1) \geq n_G - \gamma(G) + 1$.

The complete graph $G = K_n$ shows the upper bound is sharp since $n_G = n$, $\gamma(K_n) = 1$, $K_n\vee K_1\cong K_{n+1}$, and $\ZIR(K_{n+1}) = n$. By Corollary \ref{c:HveeK2},  the graph $C_r\vee K_2\cong (C_r\vee K_1)\vee K_1$ has $\ZIR(C_r\vee K_2) = r = (r+1)-1 = |V(C_r\vee K_1)| - \gamma(C_r\vee K_1)$.
\end{proof}

The \emph{wheel} graph of order $r+1$ is $W_{r+1}=C_r\vee K_1$.  

\begin{prop}\label{p:wheel}
Let $r\ge 5$. Then $\ZIR(W_{r+1})=r-\gamma(C_r)=r-\lc\frac r 3\rc$ and $\zir(W_{r+1}) =\Z(W_{r+1})=\zbar(W_{r+1})=3$.
\end{prop}
\bpf  
To establish $\ZIR(W_{r+1}) = r-\gamma(C_r)$, it suffices to show that $\ZIR(W_{r+1})\le r-\gamma(C_r)$ by Proposition \ref{p:ZIR-GveeK1}. Let $u$ denote the vertex of $K_1$, let $\hat S$ be a maximal ZIr-set of $W_{r+1}$, and for $x\in \hat S$, let $F_x$ denote private fort of $x$ with respect to $\hat S$ in $W_{r+1}$.

Suppose first that  $u\not\in F_y$ for some  $y\in \hat S\cap V(C_r)$. Then  $F_y$ is a private fort of $y$ in $C_r$ relative to $\hat S \setminus\{u\}$. Note that every fort of $C_r$ contains at least $\lc\frac r2 \rc$ vertices. If $u\not\in\hat S$, then $|\hat S|\le r-\lc\frac r 2\rc+1 = \lf\frac r 2\rf+1$. If $u\in \hat S$, then $u\not\in F_y$ for every  $y\in \hat S\cap C_r$. Thus $\hat S \setminus\{u\}$ is a ZIr-set of $C_r$ and $|\hat S|\le \ZIR(C_r)+1=\lf\frac r 2\rf+1$.  Thus in both cases $u\not\in \hat S$ and $u\in \hat S$, $|\hat S|\le\lf\frac r 2\rf+1\le  r-\lc \frac r 3\rc$ for $r\ge 5$, where the last inequality is verified algebraically.

Let $\hat S\cap V(C_r)=\{y_1,\dots,y_\ell\}$ and suppose that $u\in F_{y_i}$ for $i=1,\dots,\ell$ (which implies $u\not\in \hat S$). It suffices to show that $\ell\le r-\gamma(C_r)$. Since $u\in F_{y_i}$, $F_{y_i}$ must contain a dominating set for $C_r$ for every $y_i$. Recall that $|\hat S|\le r+1-|\cup_{i=1}^kF_{y_i}|+k=r-|\lp\cup_{i=1}^kF_{y_i}\rp\cap V(C_r)|+k$ for $k=1,\dots,\ell$ (Observation \ref{privatefortsum}).  
In order to have $|\hat S| = r-|\lp\cup_{i=1}^k F_{y_i}\rp\cap V(C_r)| + k$ there must be a set $\hat D$ such that $\hat D\cup\{y_i\}$ is a  minimum dominating set for $i=1,\dots,\ell$.  Since $r\ge 5$, there does not exist such a set for $\ell=r-\gamma(C_r)$.

It is well known that $\Z(W_{r+1})=3$ and easy to see that $\zbar(W_{r+1})=3$ since any zero forcing set must contain two consecutive cycle vertices and $u$, or three consecutive cycle vertices, and any set of this type is minimal.  By Corollary \ref{prop: min deg bnd}, $3=\delta(W_{r+1})=3\le \zir(W_{r+1})$.
\epf

We now consider the corona product of two graphs. The following notational conventions will be used throughout the rest of this section. Let $G$ and $H$ be graphs on $n_G$ and $n_H$ vertices, respectively. Let $v_1,\ldots, v_{n_G}$ denote the vertices of $G$. For $i= 1,\ldots, n_G$ let $H_i$ denote the copy of $H$ in $G\circ H$ whose vertices are adjacent to $v_i$, and let $H_i'$ denote the subgraph of $G\circ H$ induced by $V(H_i)\cup \{v_i\}$, which is isomorphic to $H\vee K_1$. The next result states basic facts about forts and ZIr-sets of coronas $G\circ H$  when $H$ has no isolated vertices. 

\begin{prop}\label{obs: forts and cut vertices}\label{obs: ZIr-sets of components} 
Let $G$ be a graph, let $H$ be a graph with no isolated vertices, and let $F$ be a fort of $G\circ H$.  Then $F'_i = F\cap V(H'_i)$ is the empty set or a fort of $H'_i$. If $S$ is a ZIr-set of $G\circ H$, then $ S'_i = S\cap V(H_i')$ is a ZIr-set of $H_i'$ and $S = \cup_{i=1}^{n_G} S'_i$.  Thus $\ZIR(G\circ H)\le n_G\ZIR(H\vee K_1)$.
\end{prop}
\bpf 
Since $H$ has no isolated vertices, every fort of $H$ has at least two elements. Thus a fort of $H_i$ is also a fort of $H'_i$. Assume  $F'_i\ne \emptyset$. If $v_i\in F'_i$, then $F'_i$ is a fort of $H'_i$.  If $v_i\not \in F'_i$, then $F'_i$ is a fort of $H_i$, so is also a fort of $H'_i$. The remaining statements follow from this property of forts.
\epf

Before proving our first result on the corona product of two graphs we introduce some new terminology. Let $S_0$ be a maximal ZIr-set of a graph $G$. An \emph{abandoned fort} relative  $S_0$  is a fort in $G$ that contains no elements of $S_0$ and 
$S_0$ \emph{abandons} a fort if there exists an abandoned fort relative to $S_0$. 
We say that $G$ \emph{abandons} a fort if there exists an upper ZIR set $S$  that abandons a fort.  

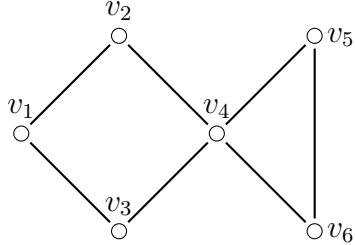
\begin{figure}[!h]
\centering
\scalebox{1}{
\begin{tikzpicture}[scale=1.3,every node/.style={draw,shape=circle,outer sep=2pt,inner sep=1pt,minimum size=.2cm}]		
\node[fill=none, label={[xshift=0pt, yshift=-5pt]$v_1$}]  (1) at (0,0) {};
\node[fill=none, label={[xshift=0pt, yshift=-5pt]$v_2$}]  (2) at (1,1) {};
\node[fill=none, label={[xshift=0pt, yshift=-5pt]$v_3$}]  (3) at (1,-1) {};
\node[fill=none, label={[xshift=0pt, yshift=-5pt]$v_4$}]  (4) at (2,0) {};
\node[fill=none, label={[xshift=10pt, yshift=-15pt]$v_5$}]  (5) at (3,1) {};
\node[fill=none, label={[xshift=10pt, yshift=-15pt]$v_6$}]  (6) at (3,-1) {};
		
\draw[thick] (4)--(3)--(1)--(2)--(4)--(5)--(6)--(4);
\end{tikzpicture}}
\caption{A graph with an upper ZIR-set that abandons a fort that is not 2-dominating.}
\label{fig: abandoned fort}
\end{figure}

\begin{ex}\label{ex:aband-fort}
Many abandoned forts naturally occur as a result of 2-dominating sets. This is not surprising since a fort is a 2-dominating set of its closed neighborhood.  Examples we have already considered include  $P_n, n\ge 3$ and $ C_n, n\ge 4$. The two vertices not in $H$ are an abandoned fort relative to $S=V(H)$ for $H\vee 2K_1$, or $H \vee K_2$ when $H \not\cong K_r$ (cf. Corollary \ref{c:HveeK2}). Figure~\ref{fig: abandoned fort} provides an example of a connected graph that abandons a fort that is not a 2-dominating set:  the upper ZIR set $\{v_1,v_4,v_5\}$ abandons the fort $\{v_2, v_3\}$, which is not 2-dominating.
\end{ex}

A single vertex forms a fort if and only if that vertex is isolated. Thus, we have the following observation.

\begin{obs}\label{obs:concede inplies order 2} 
If a maximal ZIr-set  abandons a fort $F$, then $|F|\geq 2$.
\end{obs}

\begin{thm}\label{t:HveeK1abandon}
Let $G$ and $H$ be graphs on $n_G$ and $n_H$ vertices, respectively. Assume $H$ has no isolated vertices and $H\vee K_1$ abandons a fort. Then $\ZIR(G\circ H) = n_G\ZIR(H\vee K_1)$.
\end{thm}
\begin{proof} 
By Remark~\ref{obs: ZIr-sets of components}, $\ZIR(G\circ H) \leq n_G\ZIR(H\vee K_1)$. Equality is obtained by constructing an appropriately sized ZIr-set. 

Let $S = S_1\cup \cdots \cup S_{n_G}$, where each $S_i$ is an upper ZIR set of $H_i'$ that abandons a fort $F_i$. 
Without loss of generality assume that either $v_i\in F_i$ for every $i=1,\ldots,n_G$, or $v_i\notin F_i$ for every $i=1,\ldots,n_G$. Let $u\in S_k$ for some $k  \in\{ 1,\ldots, n_G\}$  and let $F$ be a private fort of $u$ relative to $S_k$ in $H_k'$. It suffices to show $u$ has a private fort relative to $S$ in $G\circ H$. 

First, suppose $v_k\in F$. If $v_i\in F_i$ for every $i=1,\ldots,s$, then the union of $F$ with every $F_i\not= F_k$ is a private fort of $u$ relative to $S$ in $G\circ H$. If $v_i\notin F_i$ for every $i=1,\ldots,s$, then the union of $F$ with each $F_i$ such that $v_i\in N(v_k)$ is a private fort of $u$ relative to $S$ in $G\circ H$. 

Now suppose $v_k\notin F$. Then every vertex in $V(G\circ H)\setminus V(H_k')$ is adjacent to no vertices in $F$. Thus, $F$ is a private fort of $u$ relative to $S$ in $G\circ H$. Thus $S$ is a ZIr-set of $G\circ H$ and so $\ZIR(G\circ H) \geq |S| = n_G\ZIR(H\vee K_1)$.
\end{proof}

\begin{rem}\label{abandon-zfs}
Let $G$ be a graph and let $S$ be a maximal  ZIr-set set of $G$. Then $S$ abandons a fort if and only if $S$ is not a zero forcing set, because $S$ is a zero forcing set if and only if it intersects every fort by Theorem \ref{t: forts}.  Furthermore, if $G$ does not abandon a fort, then every upper ZIR set is a zero forcing set, so $\ZIR(G)=\zbar(G)$. Thus Theorem \ref{t:HveeK1abandon} applies to any graph $H$  with no isolated vertices such that $\zbar(H\vee K_1)<\ZIR(H\vee K_1)$.
\end{rem}

The next two examples provide families of graphs $H$ to which the previous theorem applies.

\begin{ex}\label{GcircC}
If $\ZIR(H\vee K_1)=n_H-\gamma(H)$ and $H\not\cong K_r$, then $H\vee K_1$ abandons a fort: Let $D_H$ be a dominating set of $H$ and let $u$ be the vertex of $K_1$.  Then the 2-dominating set $F=D_H\cup\{u\}$ is an abandoned fort of the upper ZIR set $S = V(H)\setminus D_H$ for $H\vee K_1$. In particular,  $W_{r+1}=C_r\vee K_1$ and  $\ZIR(W_{r+1}) = |V(C_r)|-\gamma(C_r)$ by Proposition \ref{p:wheel}. Thus $\ZIR(G\circ C_r)=n_G(r-\lc \frac r 3\rc)$ for $r\ge 4$.
\end{ex}

\begin{ex}\label{GcircW} 
Let $\hat H$ be a graph that is not complete. Define $H=\hat H\vee K_1$. Since $\hat H$ is not complete, neither is $H\vee K_1$, so  $\ZIR(H\vee K_1)\le |V(H\vee K_1)|-2=|V(\hat H)|$. As noted in Example \ref{ex:aband-fort}, $F=V(H\vee K_1)\setminus V(\hat H)$ is an abandoned fort relative to the upper ZIR set $S=V(\hat H)$. 
\end{ex}

Next, we present an example of a graph $H$ such that $H$ does not abandon a fort but $H\vee K_1$ does abandon a fort.

\begin{ex}  Recall that $\fr(k)=K_1\vee kK_2$ denotes the $k$th friendship graph. It was shown in Proposition \ref{p:friend} that any upper ZIR set $S$ of $\fr(k)$ must contain $v_0$ (the vertex of degree $2k$. Any fort $F_{v_0}$ of $\fr(k)$ that contains $v_0$ also contains at least one vertex from each $K_2$. However, $S$ must contain a vertex from every $K_2$ or it would not be maximal. Since every fort of $\fr(k)$ contains $v_0$ or $S$ contains either 0 or 2 for each $K_2$, $\fr(k)$ does not abandon a fort. Since $\fr(k)\vee K_1 \cong 2K_2\vee K_2$, $\fr(k)\vee K_1$ abandons a fort.
\end{ex}

We now analyze $G\circ H$ when  $H\vee K_1$ does not abandon a fort. Recall that when $H$ has no isolated vertices, every fort of $H$ contains at least two vertices, which implies the next observation. 

\begin{obs}\label{obs: lift forts in Hi}
If $H$ contains no isolated vertices, then every fort of $H_i$ is a fort of $G\circ H$.
\end{obs}

\begin{obs}\label{r:Hi}
If $n_H\ge 2$, then $\ZIR(G\circ H) \geq n_G$ because $V(H_i)$ is a fort of $G\circ H$ for $i=1,\dots,n_G$. This bound is sharp since $\ZIR(C_r\circ 2K_1) =r$ by Theorem \ref{ZIR-n/3}.
\end{obs}

Our next result 
is a lower bound that applies when  assuming only that $H$ has no isolated vertices (although Theorem \ref{t:HveeK1abandon} determines $\ZIR(G\circ H)$  when $H\vee K_1$ abandons a fort).

\begin{prop}\label{indep-LB}
Let $G$ and $H$ be graphs on $n_G$ and $n_H$ vertices, respectively. If $H$ contains no isolated vertices, then $\ZIR(G\circ H) \geq \ZIR(G)\ZIR(H\vee K_1) + (n_G-\ZIR(G))\ZIR(H)$. The bound is sharp.
\end{prop}
\begin{proof}
Let $R$ be a ZIr-set of $G$ such that $|R| = \ZIR(G)$. Without loss of generality $R = \{v_1,\ldots, v_k\}$, where $k = \ZIR(G)$. For $i = 1,\ldots,k$ let $F_i$ denote a private fort of $v_i$ relative to $S$ in $G$, and let $S_i$ be an upper ZIR set of $H_i'$. 
For $i = k+1,\ldots, {n_G}$ let $S_i$ be an upper  ZIR set of $H_i$, 
and let $T_i = V(H_i)\setminus S_i$. Let $S = S_1\cup\cdots\cup S_{n_G}$ and let $u\in S$. Then $u\in S_j$ for some $j\in \{1,\dots, {n_G}\}$. 

First, assume $j \in\{ k+1,\dots, {n_G}\}$. Then $u$ has a private fort $F$ relative to $S_j$ in $H_j$. By Observation~\ref{obs: lift forts in Hi}, $F$ is a private fort of $u$ relative to $S$ in $G\circ H$.

Now assume $j \in\{1,\ldots, k\}$. Then $u$ has a private fort $F$ relative to $S_j$ in $H_j'$. If $v_j\notin F$, then $F$ is a private fort of $u$ relative to $S$ in $G\circ H$. So, suppose $v_j\in F$ and let $T$ be the union of each $T_i$ such that $v_i\in F_j$. Then $F \cup F_j \cup T$ is a private fort of $u$ relative to $S$ in $G\circ H$. Thus, $S$ is a ZIr-set and $\ZIR(G\circ H) \geq |S| = \ZIR(G)\ZIR(H\vee K_1) + (n_G-\ZIR(G))\ZIR(H)$.

When $G=K_1$, $G\circ H=H\vee K_1$ for any $H$.  Thus $\ZIR(G)\ZIR(H\vee K_1) + (n_G-\ZIR(G))\ZIR(H)=(1)\ZIR(H\vee  K_1)+(1-1)\ZIR(H)=\ZIR(H\vee K_1)=\ZIR(G\circ H)$. 
\end{proof}

Although Proposition \ref{indep-LB} excluded isolated vertices for $H$, it is easy to compute $\ZIR(G\circ K_1)$, as in the next example.

\begin{ex}\label{Gcirc1}
Let $G$ be a graph of order $r$.  We can see that the set $S=V(G)$ is a ZIr-set of cardinality $r$: Let $u_i\in V(H_i)$ be the unique neighbor of $v_i$ that is not in $V(G)$. Then $F_i = \{v_i\}\cup \{u_j: v_j\in N[v_i]\}$ is a private fort of $v_i$. By Proposition \ref{p:2dom}  $\ZIR(G\circ K_1)\le 2r-\gamma(G\circ K_1) = 2r-r=r$, so $\ZIR(G\circ K_1)=r$.  It is known (and easy to see, e.g., by considering consider $C_{4k}$ vs. $K_r$) that $\Z(G\circ K_1)$ depends on $G$.
\end{ex}

The next result may be useful in the case where $H\vee K_1$ does not abandon a fort and $H$ has one or more isolated vertices. 

\begin{prop}
Let $G$ be a graph on $n_G$ vertices and let $H$ be a graph. Then $\ZIR(G\circ H) \geq \alpha(G)\ZIR(H\vee K_1) + (n_G-\alpha(G))(\ZIR(H) - 1)$ and this bound is sharp.
\end{prop}
\begin{proof}
Let $A$ be an independent set of $G$ such that $|A| = \alpha(G)$.  For each $v_i\in A$, let $S'_i$ be an upper ZIR set of $H_i'$. For each $v_i\in V(G)\setminus A$, let  $R_i$ be an upper ZIR set of $H_i$.  Let $u_i\in R_i$ with private fort $F_{i}$ relative to $R_i$ in $H_i$. Then $S_i = R_i\setminus\{u_i\}$ abandons the fort $F_{i}$ in $H_i$. 

Let $S = (\bigcup_{i: v_i\in A} S_i')\cup(\bigcup_{i: v_i\notin A} S_i)$. 
We show that $S$ is a ZIr-set of $G\circ H$. Let $u\in S$. Suppose first that $u\in {S_j}$ for some $j=1,\ldots,n_G$. Then $u$ has a private fort $F$ relative to $S_j$ in $H_j$ and hence, $F\cup F_j\subseteq V(H_j)$ is a private fort of $u$ relative to $S_j$ in $H_j'$ because $|F\cup F_j|\geq 2$. Now suppose $u\in S'_j$ for some $j=1,\ldots,n_G$. If $u$ has a private fort $F\subseteq V(H_j)$ relative to $S\cap V(H_j')$ in $H_j'$, then $F$ is a private fort of $u$ relative to $S$ in $G\circ H$. So, suppose this is not the case. Then $u$ has a private fort $F$ relative to $S\cap V(H_j')$ that contains $v_j$. Then the union of $F$ with every $F_i$ such that $v_i\in N(v_j)$ is a private fort of $u$ relative to $S$ in $G\circ H$. Thus, $S$ is a ZIr-set and $\ZIR(G\circ H) \geq \alpha(G)\ZIR(H\vee K_1) + (n_G-\alpha(G)) (\ZIR(H)-1).$

For sharpness, consider any graph $G$ of order $n_G$ and $H = sK_1$ with $s >n_G$. Then $\ZIR(H)=s$ and $\ZIR(H\vee K_1)=\ZIR(K_{1,s})=s-1$. It follows that $\ZIR(G\circ H) \geq \alpha(G)(s-1)+(n_G-\alpha(G))(s-1)=n_G\ZIR(H\vee K_1)=n_G(s-1)$.  Let $S$ be an upper ZIR set of $G\circ H$. Note that if $v_i\in S$, then $S\cap V(H_i')=\{v_i\}$ and therefore $|S\cap V(H_i')|\leq s$. If there exist an $H_i'$ such that $|S\cap V(H_i')|= s$, then $V(H_i)\subseteq S$ and $v_i$ is in the private fort $F_v$ of every $v\in S\cap H_i$. To ensure that $v_i$ is adjacent to two vertices of $F_v$, there exist a neighbour of $v_i$ in $G$ that is in $F_v$, say $v_k$. Since $v_k\in F_v$ it follows that $V(H_k)\subseteq F_v$. Therefore $|S|\leq (n_G-1)s=sn_G-s$ and since $s>n_G$ it contradicts the lower bound on $\ZIR(G\circ H)$. It follows that $|S\cap V(H_i)|\leq s-1$ for each $i$ and therefore $\ZIR(G\circ H) =n_G(s-1)$.
\end{proof}

\section{ Nonrelationships} 
\label{s:noncompare}
In this section we present examples that show that some parameters that are lower bounds for zero forcing number are noncomparable to the lower  ZIr number. We also show that the independence number is noncomparable to $\zir$ and $\ZIR$; this is not surprising since it is well known that it is not comparable to $\Z$, but worth noting since  $\alpha$ appears in a lower bound for $\ZIR$.  

One of the origins of zero forcing of a graph $G$ was the study of maximum nullity among the set of symmetric matrices that have graph $G$.  Let $G$ be a graph with $V(G)=\{v_1,\dots,v_n\}$. The \emph{graph} of a real symmetric $n\times n$ matrix $A$, denoted by $\G(A)$, has  vertex set $\{v_1,\dots,v_n\}$ and has edge set $\{ v_iv_j: i\neq j\mbox{ and }a_{ij} \neq 0 \}$. The \emph{maximum nullity} of the graph $G$, denoted by $\M(G)$, is defined to be \[\M(G)=\max\{\nul A:\G(A)=G\mbox{ and }A^T=A\}\] where $A^T$ denotes the transpose of $A$. It is known that $\M(G)\le \Z(G)\le \ZIR(G)$ for all $G$, so it seems natural to ask whether there is a relationship between maximum nullity and lower ZIR number. The next examples show there is no relationship between $\M$ and $\zir$.

\begin{ex}
The star $K_{1,p}$ has  $\zir(K_{1,p})=1<p-1=\M(K_{1,p})$ by Example \ref{p:star} (maximum nullity of the star is well known \cite[Theorem 9.5]{HLS22}).
\end{ex}

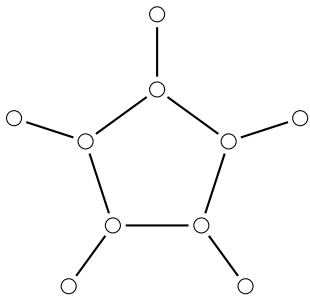
\begin{figure}[!h]
\centering
\scalebox{1}{
\begin{tikzpicture}[scale=1,every node/.style={draw,shape=circle,outer sep=2pt,inner sep=1pt,minimum size=.2cm}]		
\node[fill=none]  (00) at (0,1) {};
\node[fill=none]  (01) at (-0.951,0.309) {};
\node[fill=none]  (02) at (0.-0.588,-0.809) {};
\node[fill=none]  (03) at (0.588,-0.809) {};
\node[fill=none]  (04) at (0.951,0.309) {};
\node[fill=none]  (05) at (0,2) {};
\node[fill=none]  (06) at (-1.902,0.618) {};
\node[fill=none]  (07) at (-1.1755,-1.618) {};
\node[fill=none]  (08) at (1.1755,-1.618) {};
\node[fill=none]  (09) at (1.902,0.618) {};
		
\draw[thick] (05)--(00)--(01)--(06);
\draw[thick] (01)--(02)--(07);
\draw[thick] (02)--(03)--(08);
\draw[thick] (03)--(04)--(09);
\draw[thick] (00)--(04);
\end{tikzpicture}}
\caption{A graph $G$ having $M(G) < \zir(G)$}
\label{fig:zir>M}
\end{figure}

\begin{ex}
Let $G$ be the pentasun shown in Figure \ref{fig:zir>M}. Then it can be verified computationally that  $\zir(G)=3$ \cite{sage} and it is known that $\M(G)=2$ \cite[Example 9.118]{HLS22}.
\end{ex}

For any positive integer $k$, a \emph{$k$-tree} is $K_{k+1}$ or a graph that is obtained from $K_{k+1}$ by repeatedly adding a new vertex and joining it to an existing $k$-clique. For a graph $G$, the minimum $k$ such that $G$ is a subgraph of some $k$-tree is called the \emph{tree-width} of $G$, denoted by $\tw(G)$. It is known that $\tw(G)\le \Z(G)$ for all $G$ \cite{HLS22}, so it seems natural to ask whether there is a relationship between tree-width and lower ZIR number.  The next examples show tree-width and lower ZIR number are not comparable.

\begin{ex} 
As shown in Proposition \ref{p:friend}, the friendship graph $\fr(k)$  has $\zir(\fr(k))=k+1$. To see that $\tw(\fr(k))=2$, construct $\fr(k)$  as a subgraph of a 2-tree $G$: First define a sequence of 2-trees $L_i$ as follows: $L_1=K_3$ with $V(L_1)=\{0,1,2\}$.  For 
$i=1,\dots,  2k-2$, construct  $L_{i+1}$ 
by adding vertex $i+2$ adjacent to exactly vertices $0$ and $i+1$  of $L_i$. Observe that $V(L_i)=\{0,1,\dots,i+1\}$ and $L_i$ is a $2$ tree. From $L_{2k-1}$, delete the edges  $\{\{2i,2i+1\}:i=1,\dots,k-1\}$ to obtain $\fr(k)$.
\end{ex}

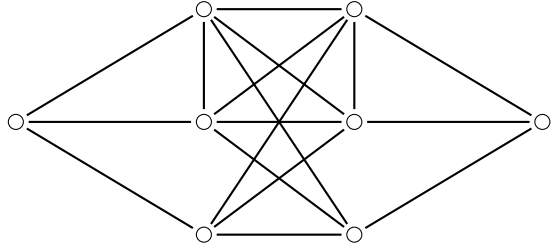
\begin{figure}[!h]
\centering
\scalebox{1}{
\begin{tikzpicture}[scale=2,every node/.style={draw,shape=circle,outer sep=2pt,inner sep=1pt,minimum size=.2cm}]		
\node[fill=none]  (00) at (-1.75,0) {};
\node[fill=none]  (01) at (-0.5,0.75) {};
\node[fill=none]  (02) at (0.5,0.75) {};
\node[fill=none]  (03) at (1.75,0) {};
\node[fill=none]  (04) at (-0.5,0) {};
\node[fill=none]  (05) at (0.5,0) {};
\node[fill=none]  (06) at (-0.5,-0.75) {};
\node[fill=none]  (07) at (0.5,-0.75) {};
		
\draw[thick] (00)--(01)--(02)--(03)--(05)--(04)--(00);
\draw[thick] (00)--(06)--(07)--(03);
\draw[thick] (01)--(04)--(02)--(05)--(01);
\draw[thick] (02)--(06)--(05);
\draw[thick] (01)--(07)--(04);
\end{tikzpicture}}
\caption{A graph $G$ having $\tw(G) > \zir(G)$}
\label{fig:tw>zir}
\end{figure}

\begin{ex}
Let $G$ be the graph shown in Figure \ref{fig:tw>zir}.  Then it has been verified computationally that $\zir(G)=4$ and  $\tw(G)=5$. \cite{sage}.
\end{ex}

Since $\pd(G)\le \Z(G)$, it is immediate that $\pd(G)\le\ZIR(G)$. Note that $C_r\circ 2K_1$ is an extremal example for power domination (least possible value of $\pd(G)$ for graphs with $3r$ vertices) and may be an extremal example for $\ZIR(G)$ (may be least possible value of $\ZIR(G)$ for graphs with $3r$ vertices). 

The next example shows that not every maximal ZIr-set is a power dominating set, but does not show that $\zir(G)<\pd(G)$ is possible.

\begin{figure}[!h]
\centering
\scalebox{1}{
\begin{tikzpicture}[scale=1.3,every node/.style={draw,shape=circle,outer sep=2pt,inner sep=1pt,minimum size=.2cm}]		
\node[fill=none, label={[xshift=-10pt, yshift=-15pt]$v_1$}]  (1) at (0.2929,0.7071) {};
\node[fill=none, label={[xshift=-10pt, yshift=-15pt]$v_2$}]  (2) at (0.2929,-0.7071) {};
\node[fill=none, label={[xshift=0pt, yshift=-5pt]$v_3$}]  (3) at (1,0) {};
\node[fill=none, label={[xshift=0pt, yshift=-5pt]$v_4$}]  (4) at (2,0) {};
\node[fill=none, label={[xshift=0pt, yshift=-5pt]$v_5$}]  (5) at (3,0) {};
\node[fill=none, label={[xshift=10pt, yshift=-15pt]$v_6$}]  (6) at (3.7071,0.7071) {};
\node[fill=none, label={[xshift=10pt, yshift=-15pt]$v_7$}]  (7) at (3.7071,-0.7071) {};

\draw[thick] (1)--(3)--(2);
\draw[thick] (3)--(4)--(5);
\draw[thick] (6)--(5)--(7);
\end{tikzpicture}}
\caption{A graph with a lower zir set that is not  power dominating.}
\label{fig: gamma and zir}
\end{figure}
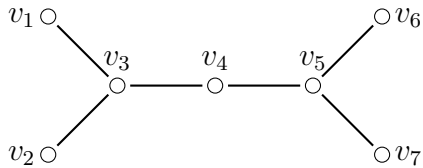

\begin{ex}
The graph $G$ in Figure \ref{fig: gamma and zir} satisfies $\zir(G) = 2$ and $\gamma_P(G) = 2$. Observe that $S = \{v_3, v_4\}$ is a maximal ZIr-set of $G$. However, $S$ is not a power dominating set.
\end{ex}

\begin{quest}\label{pd-gamma}
Is $\pd(G)\le \zir(G)$ for all graphs $G$?
\end{quest}

Well-known examples that show 
$\alpha$  is noncomparable with $\Z$  also show that  
$\alpha$  is noncomparable with $\zir$ and $\ZIR$.

\begin{ex} Consider $K_n$ with $n\ge 3$:  
$\alpha(K_n)=1 < n-1= \zir(K_n)$.
\end{ex}

\begin{ex} Consider $C_r\circ 2K_1$ with $r\ge 3$:  $\alpha(C_r\circ 2K_1)=2r>r= \ZIR(C_r\circ 2K_1)$.
\end{ex}

The domination number has appeared in several bounds for $\ZIR$. It is known that $\gamma$ is noncomparable with $\Z$. 
Consider, for example, $K_n, n\ge 3$ ($\gamma(K_n)=1$ and $\Z(K_n)=n-1$) and $C_r\circ K_1$ where $r=4\ell$ ($\gamma(C_r\circ K_1)=r$ and $\Z(C_r\circ K_1)=\frac r 2$).  What can be said $\gamma$ and $\ZIR$? From Corollary~\ref{c:ZIR2n2} a graph $G$ with $\delta(G)\geq 3$ has $\ZIR(G)\geq \frac{n}{2}$. Since a graph $G$ with no isolated vertices has $\gamma(G)\leq \frac{n}{2}$, it follows that $\gamma(G)\leq \ZIR(G)$ for graphs with minimum degree at least 3.

\begin{quest}\label{dom-gamma}
Is $\gamma(G)\le\ZIR(G)$ for all graphs $G$?
\end{quest}


\section{Computational considerations}\label{s:compute}

We have not made a careful study of the computational complexity of computing zero forcing irredundance numbers, but hope that the information in this section may help in future studies of either the complexity of computing these parameters, or in computational experimentation.

Determining all minimal forts is computationally challenging because there are examples of graph families 
such that the numbers of minimal forts is exponential in the order \cite{BMH21}. For comparison,   Genesen, Haas, and Hogben  presented a family of  trees $G_n$ of order $3n$ such that $G_n$ has $2^n$ minimum zero forcing sets and also showed that the set of minimum  zero forcing sets of a graph of order $n$ can be computed in $2^{O(n)}$ operations \cite{GHH}.

The software \cite{sage} can compute $\ZIR(G)$ and $\zir(G)$ in less than a minute for graphs up to order fifteen using \emph{Sage} installed on a 2024 MacBook Pro with an M3 Max chip and 36GB of memory running Sequoia 15.5; well-known graph families such as paths, cycles, complete graphs, and complete bipartite graphs tend to run significantly faster than random graphs.  We have not attempted to optimize the code or make it as efficient as possible.  Rather our focus was on accurately (albeit slowly) determining $\ZIR(G)$ and $\zir(G)$ for a given small graph $G$ to help in the initial investigation of examples.  

 Programs for computing zero forcing number by brute force have existed for twenty years, and over time a variety of methods for computing the zero forcing number have been developed that are significantly faster than brute force.  Such methods include the wavefront algorithm documented in Butler
et al.~\cite{wavefront} and the use of integer programming described in  \cite{BFH19}. Integer programs and their relaxations to linear programs  have also been used recently in conjunction with forts in \cite{fractZF}.  There are likely methods better than brute force for determining all forts and zero forcing irredundance numbers.


\section{Concluding remarks}\label{s:conclude} 
In this section we summarize the values of lower ZIr number, zero forcing number, upper zero forcing number, and upper ZIr number for various graph families in Table \ref{table1} and present some open questions for future work. 
In most cases listed in Table \ref{table1}, the zero forcing number was already known and most can be found (with references) in \cite[Theorem 9.5]{HLS22}.
\begin{table}[h!]
\renewcommand{\arraystretch}{1.2}
\begin{center} 
\noindent { \scriptsize\begin{tabular}{| l | c | c| c| c| c| c|}
\hline
result \#&   graph $G$ & order & $\zir(G)$ & $\Z(G)$ & $\zbar(G)$ & $\ZIR(G)$\\
\hline
\ref{e:empty} & $\OL{K_n}$  & $n$& $n$ & $n$ & $n$ & $n$ \\  
\hline
\ref{p:clique} & $K_n$  & $n$& $n-1$ & $ n-1 $ & $n-1$ & $n-1$ \\
\hline
\ref{ex:comp-pipart} & $K_{q,p},\, 1\le q \le p$  & $q+p$& $q$ & $q+p-2$  & $q+p-2$ & $q+p-2$ \\
\hline
\ref{p:path} & $P_n$, \, $5\le n$  & $n$& $1$ & $1 $ & $2$ & $\lf\frac {n-1}2\rf$\\
\hline
\ref{prop: cycle_zir_numbers} & $C_{n}, \ 4\le n$   & $n$& $2$ & $2 $ & $2$ & $\lf\frac n2\rf$ \\
\hline
\ref{p:friend} & $\fr(k), \ k\ge 2$ & $2k+1$ & $k+1$ & $k+1$ & $k+1$ & $k+1$\\
\hline
\ref{ex:all-diff} & $H(r,s), \, r\ge 2, s\ge 5, s$ odd & $r+s+1$ & 2 & $r$ & $r+1$ & $r+\frac{s-1}2$ \\
 & $H(r,3), \, r\ge 2$ & $r+s+1$ & 2 & $r$ & $r$ & $r+1$ \\
\hline
\ref{ZIR-n/2}, \ref{r:neck} & $N_k$ & $4k$ & & $k+2$ & $k+2$ & $2k$ \\
\hline
\ref{ZIR=2n/5} & $H_k$ & $5k$ & & $k+2$ &  & $2k$ \\
\hline
\ref{ZIR-n/3} & $C_r\circ 2K_1$ &3$r$ & & $r$ & $r$ & $r$ \\
\hline\ref{rK2vee2K1} & $rK_2\vee 2K_1$ & $2r+2$ & & $r+2$  &  & $2r$ \\  
\hline
\ref{prop: Join of G and H} & $ P_{n_1}\vee P_{n_2}, \ n_1,n_2\ge 7$ & $n_1+n_2$ &  & &  & $ n_1+n_2-4$ \\
\hline
\ref{p:wheel} & $W_{r+1}=C_r\vee K_1, \, r\ge 5$ & $r+1$ & $3$ & $3$ & $3$ & $r-\lc\frac r 3\rc$ \\\hline\ref{c:HveeK2} & $H\vee 2K_1$ &$|V(H)|+2$ & &  &  & $|V(H)|$ \\  
\hline
\ref{c:HveeK2} & $H\vee K_2,\, H\ne K_{|V(H)|}$ for $\ZIR$, &$|V(H)|+2$ & & $\Z(H)+2$  &   & $|V(H)|$ \\  
& $\delta(H)\ge 1$ for $\Z$ & & & & & \\
\hline
\ref{GcircC} & $H\circ C_r, \, r\ge 5$ & $|V(H)|(r+1)$ & &  &  & $ |V(H)|\lp r-\lc\frac r 3\rc\rp $ \\  
\hline
\ref{GcircW} & $H\circ W_{r+1}, \, r\ge 5$ & $|V(H)|(r+2)$ & &  &  & $ |V(H)|r $ \\  
\hline
\ref{Gcirc1} & $H\circ K_1$ & $2|V(H)|$ & &  &  & $ |V(H)| $ \\  
\hline
\end{tabular}}
\caption{Summary of values of $\zir(G), \Z(G), \zbar(G),$  $\ZIR(G)$ for various graph families $G$\vspace{-15pt}}\label{table1}
\end{center}
\end{table}

Next, we list  several  questions that were implicit in earlier discussions.

\begin{quest}
    What is the least value of $\ZIR(G)$  over all graphs of order $n$? 
\end{quest}

\begin{quest}
Theorem \ref{ZIR-n/3} provides a family of graphs $G$ satisfying $\ZIR(G)=\frac 1 3 |V(G)|$, but these graphs have $\delta(G)=1$. 
 What is an asymptotically tight lower bound on $ \ZIR(G)$ for connected graphs $G$ of  order $n$ with $\delta(G)=2$ as $n\to\infty$? 
\end{quest}
\begin{quest}
If $G$ is a connected cubic graph of order $n\ge 4$, then  $ \ZIR(G)\le \frac 3 4 n$  by Corollary \ref{c:Delta-bds}, and equality is realized by $K_4$. Do there exist connected cubic graphs $G$ of arbitrarily large order $n$ such that $ \ZIR(G)= \frac 3 4 n$?  If not, what is an asymptotically tight upper bound on $ \ZIR(G)$ for connected cubic graphs $G$ of  order $n$ as $n\to\infty$?
\end{quest}

Section \ref{s:graph-ops}  
only scratches the surface of how graph operations affect ZIr numbers. There are many related questions remaining. For example, what is the effect of graph operations on the lower zir number?

\begin{quest}
If $H$ abandons a fort, is it necessarily the case that $H\vee K_1$ abandons a fort?
\end{quest}

Section \ref{s:noncompare} contains two questions about whether the lower zir number is at least the power domination number and whether the upper ZIR number is at least the domination number.  We can also look for relationships and non-relationships with additional parameters.

 Cockayne and Mynhardt \cite{CM93} determined conditions under which integers
 $1\leq k_1\leq k_2\leq k_3\leq k_4$  can be realized as $\ir(G)=k_1, \gamma(G)=k_2, \Gamma(G)=k_3$, and $\IR(G)=k_4$. One condition is that $k_1=1$ implies that $k_2=1$; this condition would however not be true for the ZIr-numbers and zero forcing numbers since $\zir(K_{1,r})=1$, but $\Z(K_{1,r})=r-1$.  We know from Proposition \ref{ex:all-diff} that  $ k_1=2,  k_2\ge k_1$, $ k_3=k_2+1$, and $k_4> k_3$ can be realized by $H(r,s)$.

\begin{quest}
For which integer values $1\leq k_1\leq k_2\leq k_3\leq k_4$ does there exist a graph such that $\zir(G)=k_1, Z(G)=k_2, \bar{Z}(G)=k_3$ and $\ZIR(G)=k_4$?
\end{quest}

Section \ref{s:compute} discusses questions about the computational complexity of zero forcing irredundance numbers and  better methods for determining these parameters.

\begin{quest}
Determine whether the problems ``Is $\ZIR(G)\le k$?" and ``Is $\zir(G)\le k$?" are NP-complete. 
\end{quest}

\begin{quest}
Develop more efficient methods for computing $\ZIR$ and $\zir$. 
\end{quest}

\section*{Acknowledgements}

The research of B.~Curtis  is partially supported by NSF grant  1839918.  The research of A.~Roux is supported in part by the National Research Foundation of South Africa (Grant Number: 121931). The authors thank the referees for their helpful comments.


\end{document}